\documentclass{article}
\usepackage[utf8]{inputenc}
\usepackage[english]{babel}
\usepackage{amsthm}
\usepackage{graphicx}
\usepackage{placeins}
\usepackage{tikz}
\usepackage{subcaption}
\usepackage{amsmath,amssymb,xcolor,latexsym,stmaryrd}
\usepackage[hyphens]{url}
\usepackage{a4wide}
\usepackage{doi}
\usepackage{pgfplots}
\usepackage{hyperref}

\newcommand{\Th}{\mathcal{T}_{h}}
\DeclareMathOperator{\Div}{div}
\newcommand{\mbf}[1]{\boldsymbol{#1}}

\usepackage{tikz}
\usetikzlibrary{patterns}
\usepackage{pgfplots}
\pgfplotsset{compat=newest}

\newcommand{\logLogSlopeTriangle}[5]
{

    \pgfplotsextra
    {
        \pgfkeysgetvalue{/pgfplots/xmin}{\xmin}
        \pgfkeysgetvalue{/pgfplots/xmax}{\xmax}
        \pgfkeysgetvalue{/pgfplots/ymin}{\ymin}
        \pgfkeysgetvalue{/pgfplots/ymax}{\ymax}

        \pgfmathsetmacro{\xArel}{#1}
        \pgfmathsetmacro{\yArel}{#3}
        \pgfmathsetmacro{\xBrel}{#1-#2}
        \pgfmathsetmacro{\yBrel}{\yArel}
        \pgfmathsetmacro{\xCrel}{\xArel}

        \pgfmathsetmacro{\lnxB}{\xmin*(1-(#1-#2))+\xmax*(#1-#2)} 
        \pgfmathsetmacro{\lnxA}{\xmin*(1-#1)+\xmax*#1} 
        \pgfmathsetmacro{\lnyA}{\ymin*(1-#3)+\ymax*#3} 
        \pgfmathsetmacro{\lnyC}{\lnyA+#4*(\lnxA-\lnxB)}
        \pgfmathsetmacro{\yCrel}{\lnyC-\ymin)/(\ymax-\ymin)} 
        
        \coordinate (A) at (rel axis cs:\xArel,\yArel);
        \coordinate (B) at (rel axis cs:\xBrel,\yBrel);
        \coordinate (C) at (rel axis cs:\xCrel,\yCrel);

        \draw[#5]   (A)-- node[pos=0.5,anchor=north] {1}
                    (B)-- 
                    (C)-- node[pos=0.5,anchor=west] {#4}
                    cycle;
    }
}

\title{$\phi$-FEM: an efficient simulation tool using simple meshes for problems in structure mechanics and heat transfer}

\author{Stéphane Cotin$^*$, Michel Duprez\footnote{MIMESIS team, Inria Nancy - Grand Est, MLMS team, Universit\'e de Strasbourg, France.
}, Vanessa Lleras\footnote{IMAG, Univ Montpellier, CNRS, Montpellier, France.
},\\ Alexei Lozinski\footnote{Laboratoire de Math\'ematiques de Besan\c{c}on, UMR CNRS 6623, Universit\'e de Bourgogne Franche-Comt\'e, France.} \footnote{Corresponding author:  \texttt{alexei.lozinski@univ-fcomte.fr}}, and Killian Vuillemot$^*$}

\begin{document}

\maketitle






\begin{abstract}
   One of the major issues in the computational mechanics  is to take into account the geometrical complexity. To overcome this difficulty and to avoid the expensive mesh generation, geometrically unfitted methods, i.e. the numerical methods using the simple computational meshes that do not fit the boundary of the domain, and/or the internal interfaces, have been widely developed.
   In the present work, we investigate the performances of an unfitted method called $\phi$-FEM that converges optimally and uses classical finite element spaces so that it can  be easily implemented using general FEM libraries. The main idea is to take into account the geometry thanks to a level set function describing the boundary or the interface.
   Up to now, the $\phi$-FEM approach has been proposed, tested and substantiated mathematically only in some simplest settings: Poisson equation with Dirichlet/Neumann/Robin boundary conditions. Our goal here is to demonstrate its applicability to some more sophisticated governing equations arising in the computational mechanics. We consider the linear elasticity equations accompanied by either pure Dirichlet boundary conditions or by the mixed ones (Dirichlet and Neumann boundary conditios co-existing on parts of the boundary), an interface problem (linear elasticity with material coefficients abruptly changing over an internal interface), a model of elastic structures with cracks, and finally the heat equation. In all these settings, we derive an appropriate variant of $\phi$-FEM and then illustrate it by numerical tests on manufactured solutions. We also compare the accuracy and efficiency of $\phi$-FEM with those of the standard fitted FEM on the meshes of similar size, revealing the substantial gains that can be achieved by $\phi$-FEM in both the accuracy and the computational time. 
\end{abstract}

\section{Introduction}

Taking the geometrical complexity into account is one of the major issues in the computational mechanics. Although some spectacular advances in mesh generation have been achieved in  recent years, constructing and using the meshes fitting the geometry of, for example, human organs  may still be prohibitively expensive in realistic 3D configurations. Moreover, when the geometry is changing in time or on iterations of an optimization algorithm, the mesh should be frequently adapted, either by complete remeshing (expensive) or by moving the nodes (may lead to  a degradation of the mesh quality, impacting the accuracy and the stability of computations).

Geometrically unfitted methods, i.e. the numerical methods using the computational meshes that do not fit the boundary of the domain, and/or the internal interfaces, have been widely investigated in the computational mechanics for decades. Their variants come under the name of Immersed Boundary  \cite{IBMrev} or Fictitious Domain \cite{glowinski} methods. However, these classical approaches suffer from poor accuracy because of their rudimentary (but easy to implement) treatment of the boundary conditions, cf. \cite{girault}. For example, in the case of the linear elasticity equations, these methods start by extending the displacement $\mathbf{u}$, from the physical domain $\Omega$ to a fictitious domain (typically a rectangular box) $\mathcal{O}\supset\Omega$ assuming that $\mathbf{u}$ still solves the same governing equations on $\mathcal{O}$ as on $\Omega$. This creates an artificial singularity on the boundary of $\Omega$ (a jump in the normal derivative) so that the resulting numerical approximation is, at best, $\sqrt{h}$-accurate in the energy norm with whatever finite elements (from now on, $h$ denotes the mesh size).   

\begin{figure}[t]
\centering
\includegraphics[height=35mm]{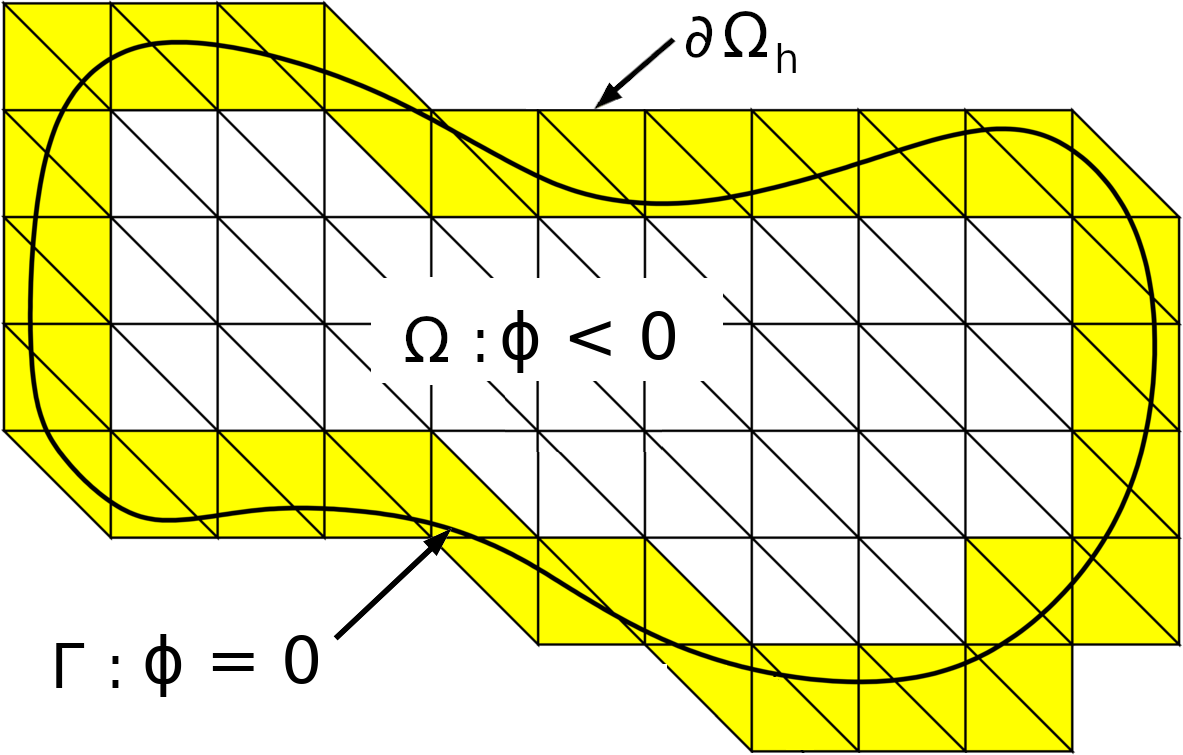}
\qquad	
\includegraphics[height=35mm]{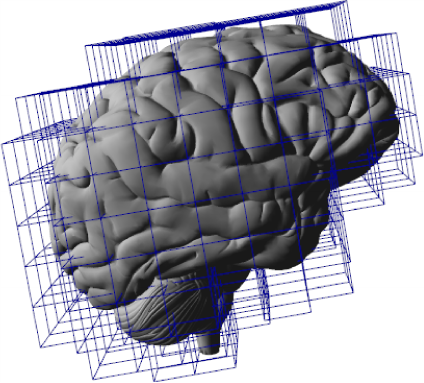}
\caption{Left: Meshes and notations for a 2D domain $\Omega=\{\phi<0\}$; the computational mesh $\Th$ is obtained from a structured background mesh and is represented by both white and yellow triangles forming the domain $\Omega_h$; the yellow triangles constitute the submesh $\Th^\Gamma$ occupying the domain $\Omega_h^\Gamma$. Right: a more involved example of an active mesh in 3D that can be used in $\phi$-FEM; a hexahedral mesh covering a brain geometry. }
\label{fig}
\end{figure}

The last two decades have seen the arrival of more accurate geometrically unfitted methods such as {XFEM} \cite{moes,haslin}, {CutFEM} \cite{burman1,burman2,cutfemrev,cutfemelasticity} and Shifted Boundary Method (SBM)  \cite{SBM,SBMelasticity}. We are citing  here only the methods based on the finite element (FE) approach; the list would be much longer if the methods based on finite differences were included. In the case of XFEM/CutFEM, the optimal accuracy, i.e. the same convergence rates as those of the standard FEM on a geometrically fitted mesh, is achieved at the price of a considerable sophistication in the  implementation of boundary conditions.  The idea is to introduce the unfitted mesh (known as the \textit{active mesh})  starting from the simple background mesh and getting rid of the cells lying entirely outside the physical domain, as illustrated at Fig. \ref{fig}. The finite elements are then set  up on the active mesh,  the variational formulation is  imposed on the \textit{physical} domain, and  an appropriate stabilization is added. In practice, one should thus compute the integrals on the actual boundary and on the parts of the active mesh cells cut by the boundary (the cut cells). To this end, one should typically construct a boundary fitted mesh, now only locally near the boundary and only for the numerical integration purposes, but the generation of a non trivial mesh is still not completely avoided. 

On the other hand, the non trivial integration is completely absent from SBM. This method introduces again an active mesh as a submesh of the background mesh (unlike CutFEM, the active mesh here contains only the cells inside $\Omega$) and then imposes the approximate boundary conditions on the boundary of the active mesh by a Taylor expansion around the actual boundary. The absence of non-standard numerical integration is an important practical advantage of SBM over XFEM/CutFEM. We note however that, to the best of our knowledge, SBM is readily available only for the lowest order FE. Moreover, in the case of Neumann boundary conditions,  the original version of SBM \cite{SBM} needs an extrapolation of the second derivatives of the solution that makes its implementation rather tricky. This difficulty can be alleviated if the problem is recast in a mixed form introducing the secondary variables for the gradient \cite{mixedSBM}. 



In this chapter, we present yet another unfitted FE-based method, first introduced in \cite{phifem,phiFEM2} and baptised $\phi$-FEM to emphasize the prominent role played in it by the level set (LS) function, traditionally denoted by $\phi$. From now on, we suppose that the physical domain is characterized by a given LS function:\footnote{In some settings presented further, the level set $\phi$ will describe an interior interface inside $\Omega$ rather than the geometry of $\Omega$ itself.}
\begin{equation}\label{Omphi}
    \Omega=\{\phi<0\}\,.
\end{equation}
Similarly to CutFEM/XFEM/SBM, we suppose that $\Omega$ is embedded into a simple background mesh and we introduce the \textit{active} computational mesh $\Th$ as in CutFEM, cf. Fig.~\ref{fig}.  However, unlike CutFEM, we abandon the variational formulation on $\Omega$. We rather introduce a non-standard formulation on the extended domain $\Omega_h$ (slightly larger than $\Omega$) occupied by the active mesh $\Th$. The general procedure is as follows:
\begin{itemize}
    \item Extend the governing equations from $\Omega$ to $\Omega_h$ and write down a formal variational formulation on $\Omega_h$ without taking into account the boundary conditions on $\partial\Omega$.
    \item Impose the boundary conditions using appropriate ansatz or additional variables, explicitly involving the level set $\phi$ which provides the link to the actual boundary. For instance, the homogeneous Dirichlet boundary conditions ($\mbf{u}=0$ on $\partial\Omega$) can be imposed by the ansatz  $\mbf{u}=\phi\mbf{w}$ thus reformulating the problem in terms of the new unknown $\mbf{w}$ (modifications for non-homogeneous conditions, mixed boundary conditions and other settings are introduced further in the text).   
    \item Add appropriate stabilization, including the ghost penalty \cite{ghost} as in CutFEM plus a least square imposition of the governing equation on the mesh cells near the boundary, to guarantee coerciveness/stability on the discrete level.
\end{itemize}
This approach allows us to achieve the optimal accuracy using classical FE spaces of any order and the usual numerical integration: all the integrals in $\phi$-FEM can be computed by standard quadrature rules on entire mesh cells and on entire boundary facets; no integration on cut cells or on the actual boundary is needed.
This is the principal advantage of $\phi$-FEM over  CutFEM/XFEM. Moreover, we can cite the following features of  $\phi$-FEM which distinguish it from both  CutFEM/XFEM and SBM:
\begin{itemize}
    \item FE of any order can be straightforwardly used in $\phi$-FEM. The geometry is naturally taken into account with the needed optimal accuracy: it suffices to approximate the LS function $\phi$ by piecewise polynomials of the same degree as that used for the primal unknown. This should be contrasted to  CutFEM  where a special additional treatment is needed if one uses FEM of order $\ge 2$. Indeed, a piecewise linear representation of the boundary is not sufficient in this case. One needs either a special implementation of the isoparametric method \cite{lehren} or a local correction by Taylor expansions \cite{boiveau}. The extension to higher order FE is not trivial for SBM either.
    \item Contrary to SBM, $\phi$-FEM is based on a purely variational formulation so that the existing standard FEM libraries suffice to implement it. The geometry of the domain comes into the formulation only through the level set $\phi$. We emphasize that $\phi$ is not necessarily the signed distance to the boundary of $\Omega$.   It is sufficient to give to the method any $\phi$ satisfying (\ref{Omphi}) which is the minimal imaginable geometrical input. This can be contrasted with SBM which assumes that the distance to the actual boundary in the normal direction is known on all the boundary facets  of the active mesh.  
\end{itemize}
Moreover, $\phi$-FEM is designed so that the matrices  of the problems on the discrete level are  reasonably conditioned, i.e. their condition numbers are of the same order as those of a standard fitting FEM on a mesh of comparable size. $\phi$-FEM shares this feature with both CutFEM/XFEM and SBM.

Up to now, the $\phi$-FEM approach has been proposed, tested and substantiated mathematically only in some simplest settings: Poisson equation with Dirichlet boundary conditions \cite{phifem}, or with Neumann/Robin boundary conditions \cite{phiFEM2}. The goal of the present chapter is to demonstrate its applicability to some more sophisticated governing equations arising in the computational mechanics. In section \ref{secElast}, we adapt $\phi$-FEM to the linear elasticity equations accompanied by either pure Dirichlet boundary conditions, or with mixed conditions (both Dirichlet and Neumann on parts of the boundary). In Section \ref{sectInterface}, we consider the interface problem (elasticity with material coefficients abruptly changing over an internal interface). Section \ref{sectFracture} is devoted to the treatment of internal cracks. Finally, our method is adapted to the heat equation in Section \ref{sectHeat}. In all these settings, we start by deriving an appropriate variant of $\phi$-FEM and then illustrate it by numerical tests on manufactured solutions. We also compare the accuracy and efficiency of $\phi$-FEM with those of the standard fitted FEM on the meshes of similar size, revealing the substantial gains that can be achieved by $\phi$-FEM in both the accuracy and the computational time.     

All the codes used in the present work have been implemented thanks to the open libraries \texttt{fenics} \cite{fenics} and \texttt{multiphenics} \cite{multiphenics}. They are available at the link\\ \url{https://github.com/michelduprez/phi-FEM-an-efficient-simulation-tool-using-simple-meshes-for-problems-in-structure-mechanics.git}

\section{Linear elasticity}\label{secElast}
In this section, we  consider the static linear elasticity for homogeneous and isotropic materials. The governing equation for the displacement $\boldsymbol u$ is thus
\begin{equation}\label{eq:elast}
\Div{\boldsymbol\sigma}(\mbf{u})+{\boldsymbol f}=0,
\end{equation}
where the stress $\boldsymbol\sigma(\mbf{u})$ is given by
$${\boldsymbol {\sigma }}(\mbf{u})=2\mu {\boldsymbol {\varepsilon }}(\mbf{u})+\lambda (\Div \mbf{u})I,
$$
$\boldsymbol{\varepsilon}(\mbf{u})=\frac{1}{2}(\nabla{\boldsymbol u} + \nabla{\boldsymbol u}^T)$ is  the strain tensor, and Lamé parameters $\lambda,\mu$ are defined via the Young modulus $E$ and the Poisson coefficient $\nu$ by
\begin{equation}\label{muEnu}
    \mu = \dfrac{E}{2(1+\nu)}\mbox{ and }\lambda=\dfrac{E\nu}{(1+\nu)(1-2\nu)}\,.
\end{equation}
Equation (\ref{eq:elast}) is posed in a domain $\Omega$, which can be two or three dimensional, and should be accompanied with Dirichlet and Neumann boundary conditions on $\Gamma=\partial\Omega$. We assume that $\Gamma$ is decomposed into two disjoint parts, $\Gamma=\Gamma_D\cup\Gamma_N$ with $\Gamma_D\neq \varnothing$, and
\begin{align} 
\mbf{u}&=\mbf{u}^g \text{ on } \Gamma_D, \label{bcGammaD} \\
\mbf{\sigma}(\mbf{u})\mbf{n}&=\mbf{g} \ \text{ on } \Gamma_N, \label{bcGammaN} 
\end{align}
with the given displacement $\mbf{u}^g$ on $\Gamma_D$ and the given force $\mbf{g}$ on $\Gamma_N$.

Let us first recall the weak formulation of this problem (to be compared with
forthcoming $\phi$-FEM formulations): find the vector field $\mbf{u}$ on $\Omega$ s.t. $\mbf{u}
|_{\Gamma_D} =\mbf{u}^g $ and
\begin{equation}\label{weakSimple}
  \int_{\Omega} \mbf{\sigma} (\mbf{u}) : \nabla\mbf{v}   
   = \int_{\Omega} \mbf{f} \cdot \mbf{v}+
   \int_{\Gamma_N} \mbf{g} \cdot \mbf{v}, \quad \forall
   \mbf{v} \text{ on }\Omega \text{ such that } \mbf{v} |_{\Gamma_D} = 0  .
\end{equation}
This is obtained by multiplying the equation by a test function
$\mbf{v}$, integrating over $\Omega$ and taking into account the boundary
conditions. Formulation (\ref{weakSimple}) is routinely used to construct conforming FE methods, which necessitate a mesh that fits the domain $\Omega$ in order to approximate the integrals on $\Omega$ and $\Gamma_N$ and to impose $\mbf{u} =\mbf{u}^g $ on ${\Gamma_D}$. 

We now consider the situation where a fitting mesh of $\Omega$ is not available. We rather assume that $\Omega$ is inscribed in a box ${\mathcal{O}}$ which is covered by a simple background mesh $\Th^{\mathcal{O}}$. We further introduce the computational mesh $\Th$ (also referred to as the active mesh) by getting rid of cells lying entirely outside $\Omega$. In practice,  $\Omega$ is given by the level-set function $\phi$: $\Omega=\{\phi<0\}$. Usually, the level set is known only approximately. Accordingly, we assume that we are given a FE function $\phi_h$, i.e. a piecewise polynomial function on mesh $\Th^{\mathcal O}$, which approximate sufficiently well $\phi$. The selection of the mesh cells forming the active mesh is done on the basis of $\phi_h$ rather than $\phi$: 
\begin{equation}\label{Th}
\Th:=\{T\in \Th^{\mathcal{O}}:T\cap \{\phi_h<0\}\neq  \varnothing\}\,.
\end{equation}
The domain occupied by $\Th$ is  denoted by $\Omega_h$, \textit{i.e.}
${\Omega_h} = (\cup_{T \in \Th}T)^o$. In some of our methods, we shall also need a submesh of $\Th$, referred to as  $\Th^\Gamma$, consisting of the cells intersected with the curve (surface) $\{\phi_h=0\}$, approximating $\Gamma$:
\begin{equation}\label{ThGam}
\Th^\Gamma:=\{T\in \Th^{\mathcal{O}}:T\cap \{\phi_h=0\}\neq  \varnothing\}\,.
\end{equation}
The domain covered by mesh $\Th^\Gamma$ will be denoted by $\Omega_h^\Gamma$, cf. Fig.~\ref{fig}. 

The starting point of all variants of $\phi$-FEM is a variational formulation of problem (\ref{eq:elast}) extended to $\Omega_h$, in which we do not impose any boundary conditions since they are lacking on $\partial\Omega_h$. We thus assume that the right-hand side $\mbf{f}$ is given on the whole $\Omega_h$ rather than on $\Omega$ alone, and suppose moreover that $\mbf{u}$ can be extended from $\Omega$ to $\Omega_h$ as  the solution to the  governing equation (\ref{eq:elast}), now posed on $\Omega_h$ instead of $\Omega$. In a usual manner, we take then any test function
$\mbf{v}$ on $\Omega_h$, multiply the governing equation by $\mbf{v}$ and
integrate it over $\Omega_h$. This gives the following formulation: find a vector field $\mbf{u}$ on $\Omega_h$ such that
\begin{equation}\label{weakExt} 
\int_{\Omega_h} \mbf{\sigma} (\mbf{u}) : \nabla\mbf{v}  
- \int_{\partial\Omega_h} \mbf{\sigma} (\mbf{u})  \mbf{n} \cdot \mbf{v}
= \int_{\Omega_h} \mbf{f} \cdot \mbf{v}, \quad \forall \mbf{v} \text{ on }\Omega_h
\end{equation}
We emphasize that this formulation is fundamentally different from the standard formulation (\ref{weakSimple}). First of all, no boundary conditions are incorporated in (\ref{weakExt}) so that we cannot expect it to admit a unique solution. Furthermore, if we add somehow the boundary conditions on $\partial\Omega$ to (\ref{weakExt}), which we shall do indeed when constructing our $\phi$-FEM variants, the resulting formulation will still be ill posed, meaning that its solution (on the continuous level) either does not exist, or is not unique. However, we shall be able to turn these problems into well defined numerical schemes by adding an appropriate stabilization on the discrete level. 

\subsection{Dirichlet conditions}
Let us first consider the case of pure Dirichlet conditions: $\Gamma=\Gamma_D$. 
On the continuous level, we want thus to impose $\mbf{u}=\mbf{u}^g$ on $\Gamma = \Gamma_D = \{ \phi = 0 \}$ on top of the general formulation (\ref{weakExt}) of the problem on $\Omega_h$. We consider here 2 options to achieve this: 1) \textbf{direct} Dirichlet $\phi$-FEM, as proposed in \cite{phifem}, introducing a new unknown $\mbf{w}$ and redefining $\mbf{u}$ through the product $\phi\mbf{w}$ which automatically vanishes on $\Gamma$; 2) \textbf{dual} Dirichlet $\phi$-FEM, inspired by \cite{phiFEM2}, keeping the original unknown $\mbf{u}$ and imposing $\mbf{u}=\mbf{u}^g$ on $\Gamma$ with the aid of an auxiliary variable $\mbf{p}$ in a least-square manner. In more details, our two approaches can be described as follows:    
\begin{itemize}
  \item \textbf{Direct Dirichlet $\phi$-FEM} \textit{(on continuous level)}. Supposing that $\mbf{u}^g$ is actually given on the whole $\Omega_h$ rather than on $\Gamma$ alone, we make the ansatz 
  \begin{equation}\label{ansatz}
  	\mbf{u}=\mbf{u}^g + \phi \mbf{w},\text{ on }\Omega_h
  \end{equation} 
  and substitute it into (\ref{weakExt}). To make the formulation more symmetric we also
  replace the test functions $\mbf{v}$ by $\phi \mbf{z}$. This yields: find a vector field $\mbf{w}$ on $\Omega_h$ such that
  \begin{multline}\label{DirectPhi} \int_{\Omega_h} \mbf{\sigma} (\phi\mbf{w}) : \nabla(\phi\mbf{z})  
  - \int_{\partial \Omega_h} \mbf{\sigma}  (\phi \mbf{w}) \mbf{n} \cdot \phi
     \mbf{z}
     = \int_{\Omega_h} \mbf{f} \cdot \phi \mbf{z} \\
     - \int_{\Omega_h} \mbf{\sigma}  (\mbf{u}^g)
     : \nabla (\phi \mbf{z})  + \int_{\partial
     \Omega_h} \mbf{\sigma}  (\mbf{u}^g) \mbf{n}
     \cdot \phi \mbf{z}, \quad 
     \forall \mbf{z}\text{ on }\Omega_h .
  \end{multline}
  The idea is thus to work with the new unknown $\mbf{w}$ on $\Omega_h$,
  discretize it by FEM starting from the variational formulation above, and to
  reconstitute the approximation to $\mbf{u}$ by the ansatz (\ref{ansatz}).
  
  \item \textbf{Dual Dirichlet $\phi$-FEM}  \textit{(on continuous level)}.     We now suppose that $\mbf{u}^g$ is  defined on $\Omega_h^{\Gamma}$, cf. (\ref{ThGam}),
  rather than on the whole of $\Omega_h$. We keep the primal unknown
  $\mbf{u}$ in (\ref{weakExt}) and we want to impose
  \begin{equation}\label{ansatzDual}
  	\mbf{u}=\mbf{u}^g + \phi \mbf{p},\text{ on }\Omega_h^\Gamma
  \end{equation} 
  on top of it, with a new auxiliary unknown $\mbf{p}$ on $\Omega_h^\Gamma$. The new variable $\mbf{p}$ lives beside $\mbf{u}$ inside a variational formulation that combines (\ref{weakExt}) with (\ref{ansatzDual}):  find vector fields $\mbf{u}$ on $\Omega_h$ and   $\mbf{p}$ on $\Omega^{\Gamma}_h$ such that
  \begin{multline}\label{DualPhi}  
    \int_{\Omega_h} \mbf{\sigma}  (\mbf{u}) : \nabla\mbf{v}   
     - \int_{\partial\Omega_h} \mbf{\sigma}  (\mbf{u}) \mbf{n} \cdot \mbf{v} 
     + {\gamma} \int_{\Omega_h^{\Gamma}} (\mbf{u}- \phi \mbf{p}) \cdot
     (\mbf{v}- \phi \mbf{q}) \\
     = \int_{\Omega_h} \mbf{f} \cdot
     \mbf{v}+  {\gamma}\int_{\Omega_h^{\Gamma}} \mbf{u}^g \cdot
     (\mbf{v}- \phi \mbf{q}), \quad \forall \mbf{v} \text{ on }
     \Omega_h, \mbf{q} \text{ on } \Omega_h^{\Gamma} 
  \end{multline} 
  with a positive parameter $\gamma$. Comparing the direct and dual variants, we observe that the expressions (\ref{ansatz}) and (\ref{ansatzDual}) are of course pretty similar, but their roles are quite different in the corresponding methods. The variable $\mbf{w}$ replaces $\mbf{u}$ in (\ref{DirectPhi}), while $\mbf{p}$ lives alongside $\mbf{u}$ in (\ref{DualPhi}). The introduction of the additional variable $\mbf{p}$  makes the dual method only slightly more expensive than the direct one, since this new variable is introduced only on a narrow strip around $\Gamma$. On the other hand, a certain advantage of the dual variant over the  direct one lies in the fact that both $\phi$ and $\mbf{u}^g$ should be here known only locally around $\Gamma$ since they enter into equation (\ref{DualPhi}) only on $\Omega_h^{\Gamma}$. This can facilitate the construction of $\phi$ and $\mbf{u}^g$ in practice. More importantly, it is the dual method that we shall be able to adapt to various, more and more complicated settings below.  
\end{itemize}
As mentioned above, both variational problems (\ref{DirectPhi}) and (\ref{DualPhi}) are derived on a very formal level. They are not valid in any mathematically rigorous way: we cannot expect to have a meaningful boundary value problems on a domain $\Omega_h$ with no
boundary conditions on $\partial \Omega_h$, while prescribing some conditions
on a curve (surface) $\Gamma$ which is inside $\Omega_h$. However, both
formulations can serve as starting problems to write down FE problems which
become well-posed once an appropriate stabilization is added.  

We start by introducing the FE spaces: fix an integer $k\ge 1$ and let
\begin{equation} \label{spaceVh}
 V_{h} := \left\lbrace \mbf{v}_h :\Omega_{h}\to\mathbb{R}^d : \mbf{v}_{h |T} \in \mathbb{P}^k(T)^d \ \ \forall T \in \mathcal{T}_h, \ \mbf{v}_h \text{ continuous on }\Omega_h \right\rbrace .
 \end{equation}
For future reference, we introduce the local version of this space for any submesh $\mathcal{M}_h$ of $\Th$ and polynomial degree $l\ge 0$
\begin{equation}
 \label{spaceQh}
 Q_{h}^l(\mathcal{M}_h) := \left\lbrace \mbf{q}_h:\mathcal{M}_{h} \to\mathbb{R}^d : \mbf{q}_{h |T} \in \mathbb{P}^{l}(T)^{d} \ \ \forall T \in \mathcal{M}_h , \ \mbf{q}_h \text{ continuous on }\mathcal{M}_h\text{ if }l\ge 0\right\rbrace.
 \end{equation}
In particular, we shall need the space   $Q_{h}^k(\Omega_{h}^{\Gamma})$ on the submesh $\Omega_{h}^{\Gamma}$ in the Dual version of Dirichlet $\phi$-FEM.
 
The two variants of $\phi$-FEM introduced above can now be written on the fully
discrete level  as: 
 \begin{itemize}
  \item \textbf{Direct Dirichlet $\phi$-FEM}: find $\mbf{w}_h\in V_h$ such that
\begin{multline}\label{DirectDiscrete} \int_{\Omega_h} \mbf{\sigma}  (\phi_h \mbf{w}_h) :
\nabla (\phi_h \mbf{z}_h) 
- \int_{\partial \Omega_h}
\mbf{\sigma}  (\phi_h \mbf{w}_h) \mbf{n} \cdot \phi_h\mbf{z}_h 
+ G_h(\phi_h\mbf{w}_h,\phi_h\mbf{z}_h)  
+ J_h^{lhs} (\phi_h\mbf{w}_h,\phi_h\mbf{z}_h) \\ 
= \int_{\Omega_h} \mbf{f} \cdot \phi_h \mbf{z}_h 
- \int_{\Omega_h} \mbf{\sigma}  (\mbf{u}_h^g)
: \nabla (\phi_h \mbf{z}_h) 
+ \int_{\partial\Omega_h} \mbf{\sigma}  (\mbf{u}_h^g) \mbf{n}
\cdot \phi_h \mbf{z}_h, \\ 
 + J_h^{rhs} (\phi_h\mbf{z}_h),\quad\forall \mbf{z}_h\in V_h 
\end{multline}
and set $\mbf{u}_h=\mbf{u}_h^g + \phi_h \mbf{w}_h$. Here $\phi_h, \mbf{u}^g_h$ are FE
approximations  for $\phi, \mbf{u}^g$ on the whole $\Omega_h$, and $G_h,J_h^{lhs}, J_h^{rhs}$ stand for the stabilization terms 
\begin{equation}\label{Gh}
G_h(\mbf{u}, \mbf{v}) : = \sigma_D h \sum_{E \in \mathcal{F}_h^{\Gamma}} \int_E \left[ \mbf{\sigma}(\mbf{u})\mbf{n}  \right] \cdot \left[ \mbf{\sigma}(\mbf{v})\mbf{n}\right] ,
\end{equation} 
\begin{equation}\label{Jhlrhs}
J_h^{lhs}(\mbf{u}, \mbf{v}) : =  \sigma_D h^2  \sum_{T \in \mathcal{T}_h^{\Gamma}} \int_T  \Div \mbf{\sigma}(\mbf{u})\cdot \Div\mbf{\sigma}(\mbf{v}) \,,\qquad
J_h^{rhs} (\mbf{v}) : =- \sigma_D h^2\sum_{T \in \Th^{\Gamma}} \int_T \mbf{f} \cdot\Div\mbf{\sigma}(\mbf{v}) 
\,.
\end{equation}
The stabilization $G_h$ (\ref{Gh}) is known as the ghost penalty. $\sigma_D$ in (\ref{Gh}) is a positive stabilization parameter which should be chosen sufficiently big (in a mesh independent manner). $\mathcal{F}_h^{\Gamma}$ stands for the set of internal facets of mesh $\Th$ which are also the facets of $\Th^\Gamma$ (these are the facets either intersected by $\Gamma$, or belonging to the cells intersected by $\Gamma$).  Stabilization (\ref{Gh}) was first introduced in \cite{ghost} in the form of penalization of jumps in the normal derivatives of the FE solution. Here,  we prefer to penalize the jumps of internal elastic forces, following \cite{claus2018stable}, thus controlling appropriate combinations of the derivatives, rather than the normal derivatives themselves.   We emphasize however that the original ghost penalty in \cite{cutfemelasticity} also involved the jumps of higher order derivatives of $\mbf{u}$ (up to the highest order of polynomials present in the FE formulation), while our variant affects the first order derivatives only. We can allow ourselves to reduce the order of stabilized derivatives thanks to the presence of additional stabilization terms $J_h^{lhs}$ (\ref{Jhlrhs}), as first suggested in \cite{phifem} (a similar idea can also be found in \cite{elfverson2019new}). The combination of $G_h$ and $J_h^{lhs}$ allows us indeed to get rid of possible spurious oscillations of the approximate solution on ``badly cut" cells near $\Gamma$ and to guarantee the coerciveness of the bilinear form in our FE formulation.  Note that the terms $J_h^{lhs}$ are not consistent by themselves but they are  consistently compensated  by their right-hand side counterpart $J_h^{rhs}$. Indeed, the exact solution satisfies  $\Div\mbf{\sigma}(\mbf{u})=-\mbf{f}$ so that $J_h^{lhs}(\mbf{u}, \mbf{v})=J_h^{rhs}(\mbf{v})$ if $\mbf{u}$ is the exact solution.

  \item Dual $\phi$-FEM-Dirichlet: find $\mbf{u}_h\in V_h$, $\mbf{p}_h\in Q_{h}^k(\Omega_{h}^{\Gamma})$  such that
  \begin{multline}  \label{DualDiscrete}
  \int_{\Omega_h} \mbf{\sigma}  (\mbf{u}_h) : \nabla\mbf{v}_h 
  - \int_{\partial\Omega_h} \mbf{\sigma}(\mbf{u}_h) \mbf{n} \cdot \mbf{v}_h 
  + \frac{\gamma}{h^2} \int_{\Omega_h^{\Gamma}} (\mbf{u}_h- \frac{1}{h}\phi_h \mbf{p}_h) \cdot
(\mbf{v}_h- \frac{1}{h}\phi_h \mbf{q}_h) 
\\+ G_h (\mbf{u}_h,\mbf{v}_h) 
+ J_h^{lhs} (\mbf{u}_h,\mbf{v}_h) \\
= \int_{\Omega_h} \mbf{f} \cdot\mbf{v}_h
+ \frac{\gamma}{h^2} \int_{\Omega_h^{\Gamma}} \mbf{u}_h^g \cdot
(\mbf{v}_h- \frac{1}{h}\phi_h \mbf{q}_h)
+ J_h^{rhs} (\mbf{v}_h), \quad \forall \mbf{v}_h \in V_h,\ 
\mbf{q}_h \in Q_{h}^k(\Omega_{h}^{\Gamma}).
\end{multline} 
With respect to (\ref{DualPhi}, we have added here the factors $\frac{1}{h}$, $\frac{1}{h^2}$. They serve to control the condition numbers, cf. \cite{phiFEM2}. The stabilizations $G_h$, $J_h^{lhs}$, $J_h^{rhs}$ are again defined by (\ref{Gh}) and (\ref{Jhlrhs}).
\end{itemize}

\begin{figure}[htbp]
\centering 
 \includegraphics[width = 0.48\textwidth]{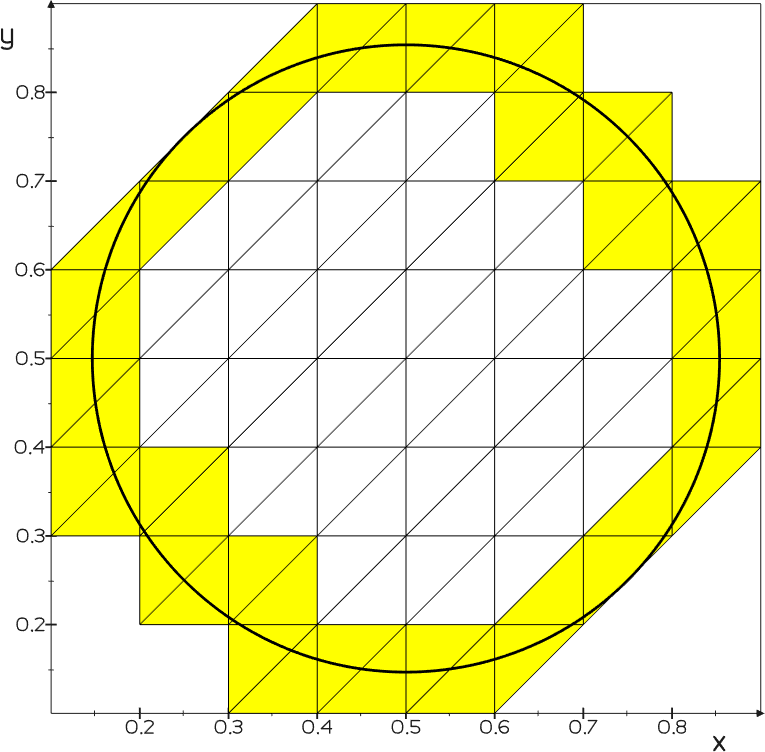}
\quad \includegraphics[width=0.48\textwidth]{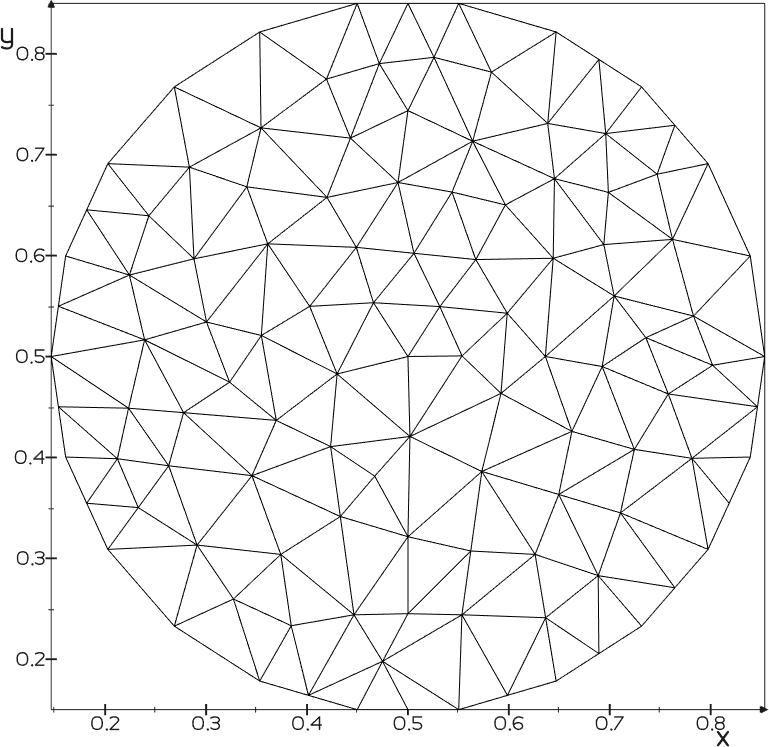}
\caption{Circular domain given by (\ref{phiCircle}).  Left: active meshes for $\phi$-FEM (with cells from $\Th^\Gamma$ in yellow). Right: a fitted mesh for the standard FEM. }\label{fig:meshes dirichlet}
\end{figure}

\begin{figure}[htbp]
\centering
\begin{tikzpicture}
\begin{loglogaxis}[name = ax1, width = .45\textwidth, xlabel = $h$, 
            ylabel = $L^2$ relative error,
            legend style = { at={(0.7,1.3)}, legend columns =1,
			/tikz/column 2/.style={column sep = 10pt}}]
\addplot coordinates {
(0.08838834764831845,1.5588115850660593e-06)
(0.04419417382415922,9.469566130085001e-08)
(0.02209708691207961,7.2400075783873776e-09)
(0.011048543456039806,8.204951465874365e-10)
(0.005524271728019903,1.0117386213743662e-10) };
\addplot coordinates { 
(0.08838834764831845,0.00017927946943448)
(0.04419417382415922,9.288747432430552e-06)
(0.02209708691207961,8.561278085915226e-07)
(0.011048543456039806,6.216451875056968e-08)
(0.005524271728019903,6.4733637784675345e-09)};
\addplot coordinates { 
(0.08814237266741118,0.0013863242202072467)
(0.044176082108213596,0.0003411362090257958)
(0.022095901417288732,8.270844843658769e-05)
(0.01104624219043868,2.0738926136251253e-05)
(0.005523930001854928,5.168354870491673e-06)};
\logLogSlopeTriangle{0.53}{0.2}{0.12}{3}{blue};
\legend{Direct $\phi$-FEM,Dual $\phi$-FEM, Standard FEM}
\end{loglogaxis}
\end{tikzpicture}
\begin{tikzpicture}
\begin{loglogaxis}[name = ax1, width = .45\textwidth, xlabel = $h$, 
ylabel = $H^1$ relative error,
legend style = { at={(0.7,1.3)}, legend columns =1,
/tikz/column 2/.style={column sep = 10pt}}]

\addplot coordinates { 
(0.08838834764831845,4.4726648824471084e-05)
(0.04419417382415922,6.698604600427322e-06)
(0.02209708691207961,1.4472949769396173e-06)
(0.011048543456039806,3.5820477725646584e-07)
(0.005524271728019903,8.9819527477158e-08)};

\addplot coordinates{
(0.08838834764831845,0.0023935969083683033)
(0.04419417382415922,0.00027798873140552567)
(0.02209708691207961,5.4548588878934805e-05)
(0.011048543456039806,9.040565425995086e-06)
(0.005524271728019903,1.8279013328542586e-06)};

\addplot coordinates{
(0.08814237266741118,0.00874593661775123)
(0.044176082108213596,0.0031309191906407226)
(0.022095901417288732,0.0011619674908133126)
(0.01104624219043868,0.0003941847087612613)
(0.005523930001854928,0.00013882795038475216)};
\logLogSlopeTriangle{0.53}{0.2}{0.12}{2}{blue};
\legend{Direct $\phi$-FEM,Dual $\phi$-FEM, Standard FEM}
\end{loglogaxis}
\end{tikzpicture}
\caption{Test case with pure Dirichlet conditions.  $L^2$  relative errors on the left, $H^1$ relative errors on the right.}\label{fig:diriclet_elasticity1}
\end{figure}
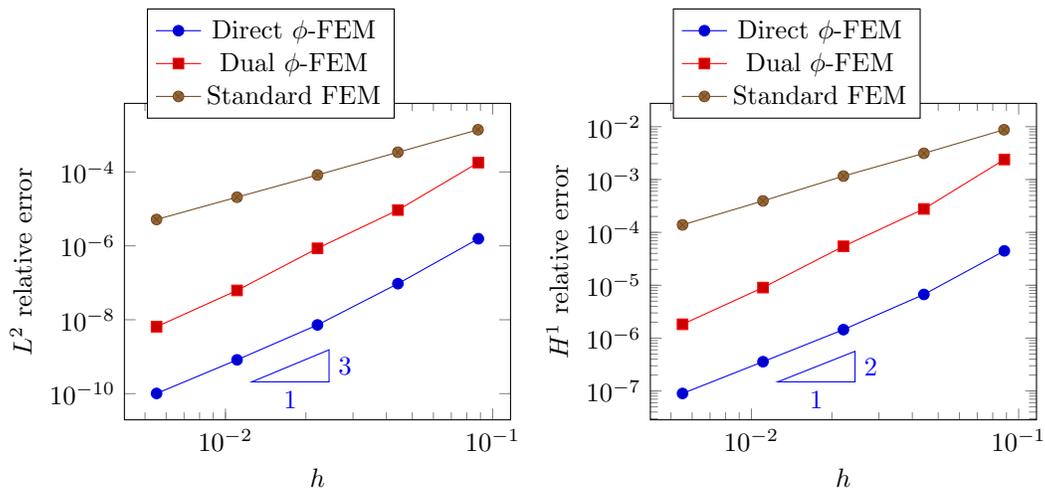

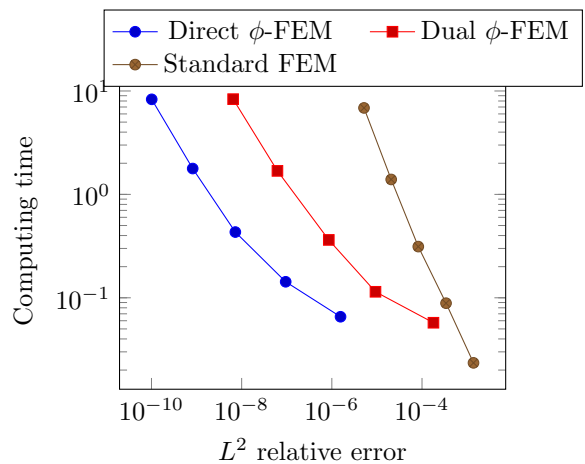
\begin{figure}[htbp]
\centering
\begin{tikzpicture}
\begin{loglogaxis}[name = ax1, width = .45\textwidth, xlabel = $L^2$ relative error, 
ylabel = Computing time,
anchor=north west,legend style = { at={(1.2,1.2)}, legend columns =2,
/tikz/column 2/.style={column sep = 10pt}}]
\addplot coordinates {
(1.5588115850660593e-06,0.06556248664855957)
(9.469566130085001e-08,0.1428232192993164)
(7.2400075783873776e-09,0.43224072456359863)
(8.204951465874365e-10,1.778151512145996)
(1.0117386213743662e-10,8.303027629852295)
};
\addplot coordinates {
(0.00017927946943448,0.0573728084564209)
(9.288747432430552e-06,0.11432981491088867)
(8.561278085915226e-07,0.36284589767456055)
(6.216451875056968e-08,1.687530279159546)
(6.4733637784675345e-09,8.336758613586426)};
\addplot coordinates { 
(0.0013863242202072467,0.023479223251342773)
(0.0003411362090257958,0.08866548538208008)
(8.270844843658769e-05,0.31266045570373535)
(2.0738926136251253e-05,1.3940625190734863)
(5.168354870491673e-06,6.85805869102478)};
\legend{Direct $\phi$-FEM, Dual $\phi$-FEM, Standard FEM}
\end{loglogaxis}
\end{tikzpicture}
\caption{Test case with pure Dirichlet conditions.  Computing time (in seconds) vs. the $L^2$  relative errors.}\label{fig:computedtime_elasticity}
\end{figure}

\paragraph*{Test case:}\label{testcaseDir}
Let $\mathcal{O}$ be the square $(0,1)^2$ and $\mathcal{T}_h^{\mathcal{O}}$ a uniform mesh on $\mathcal{O}$. Let $\Omega$ be the circle centered at the point $(0.5,0.5)$ of radius $\frac{\sqrt{2}}{4}$. The level set function $\phi$ is thus given by 
\begin{equation}\label{phiCircle}
    \phi(x,y) = -\frac{1}{8} + (x-0.5)^2 + (y-0.5)^2\,. 
\end{equation}
We take the elasticity parameters $E=2$ and $\nu=0.3$, 
and the scheme parameters $\gamma = \sigma_D = 20.0$. We use $\mathbb{P}^2$-Lagrange polynomials  for both FE spaces $V_h$ and $Q_h$, i.e. we set $k=2$ in (\ref{spaceVh}) and (\ref{spaceQh}).
We finally choose a manufactured exact solution  
\begin{equation}\label{manufu}
    \mbf{u}=\mbf{u}_{ex} := (\sin(x)  \exp(y), \sin(y)  \exp(x))
\end{equation} 
giving the right hand side $\mbf{f}$ by substitution to (\ref{eq:elast}) and the
boundary conditions $\mbf{u}^g=\mbf{u}_{ex}$ on $\Gamma$. In order to set up both $\phi$-FEM schemes above, we should extend $\mbf{u}^g$ from $\Gamma$ to $\Omega_h$ (in the case of the direct method) or to $\Omega_h^\Gamma$ (in the case of the dual method). To mimic the realistic situation where  $\mbf{u}^g$ is known on $\Gamma$ only, we prefer not to extend $\mbf{u}^g$ by $\mbf{u}_{ex}$ everywhere. We rather set
$$
\mbf{u}^g = \mbf{u}_{ex}(1+\phi), \quad \text{ on } \Omega_h\text{ or on }\Omega_h^{\Gamma}
$$
adding to  $\mbf{u}_{ex}$ a  perturbation which vanishes on the boundary.  

The typical active meshes $\Th$ and $\Th^\Gamma$ for $\phi$-FEM are illustrated on Fig.~\ref{fig:meshes dirichlet} (left).  Besides the direct $\phi$-FEM (\ref{DirectDiscrete}) and the dual $\phi$-FEM (\ref{DualDiscrete}), we shall present the numerical results obtained by the standard FEM with $\mathbb{P}^2$-Lagrange polynomials on fitted meshes for approximately the same values of $h$, as illustrated on Fig.~\ref{fig:meshes dirichlet} (right). The results obtained by both variants of $\phi$-FEM and by the standard FEM are reported in Figs.~\ref{fig:diriclet_elasticity1} and \ref{fig:computedtime_elasticity}.

We first illustrate the numerical convergences order for the relative errors in $L^2$ and $H^1$ norms at Fig.~\ref{fig:diriclet_elasticity1}. We observe that both variants of $\phi$-FEM demonstrate indeed the expected optimal convergence orders: $h^2$ is the $H^1$-seminorm and $h^3$ in the $L^2$-norm, and the direct variant performs significantly better than the dual one.  This can be attributed to a better representation of the solution near the boundary in the direct variant: indeed it is effectively approximated there by fourth-order polynomials ($\mathbb{P}^2$ for $\mbf{w}_h$ times $\mathbb{P}^2$ for $\phi_h$). Moreover, both $\phi$-FEMs, even the dual one, significantly outperform the standard FEM (the latter is even of a suboptimal order in the $L^2$-norm). This can be partially attributed to a coarse geometry approximation. Indeed, we use triangular meshes so that the curved boundary of $\Omega$ is actually approximated by a collection  of straight segments, i.e. the boundary facets of the fitted mesh, cf. Fig.~\ref{fig:meshes dirichlet} (right).  The superior efficiency of $\phi$-FEM with respect  to the standard FEM is further confirmed by Fig.~\ref{fig:computedtime_elasticity}. We report there the  computing times on different meshes for the 3 methods and set them against the relative $L^2$ error. These computing times include assembling of the FE matrices and resolution of the resulting linear systems. For a given relative error, the calculations are always much faster with $\phi$-FEM than with the standard FEM. The advantage would be even more significant if the mesh generation times were included, since the construction of active meshes in $\phi$-FEM only involves choosing a subset of cells according to a simple criterion, and some renumbering of the degrees of freedom. We do not dispose however of an efficient implementation of cell selection algorithm at the moment.  All our computations are performed using the Python interface for the popular FEniCS computing platform, and the selection of active cells is done by a simple, non-optimized Python script.

\subsection{Mixed boundary conditions}


We now consider the much more complicated case of mixed conditions (\ref{bcGammaD})--(\ref{bcGammaN}) on the boundary $\Gamma = \Gamma_N \cup \Gamma_D$ with $\Gamma_D\neq \varnothing$ and $\Gamma_N\neq \varnothing$. This setting is challenging for any geometrically unfitted method since the junction between the Dirichlet and Neummann boundary parts can occur inside a mesh cell, so that approximating polynomials in this cell should account simultaneously for both boundary conditions.  In \cite{cutfemelasticity}, it is demonstrated that the linear elasticity with mixed boundary conditions can be successfully treated by CutFEM. A rigorous mathematical substantiation allowing of the low regularity of the solution is available in \cite{burmanmixed}. 
Here, we shall adapt  $\phi$-FEM (in the dual form) to the mixed boundary conditions by adopting a ``lazy" approach: we choose to do not impose any boundary conditions on a mesh cell if the Dirichlet/Neumann junction happens to be inside it. 

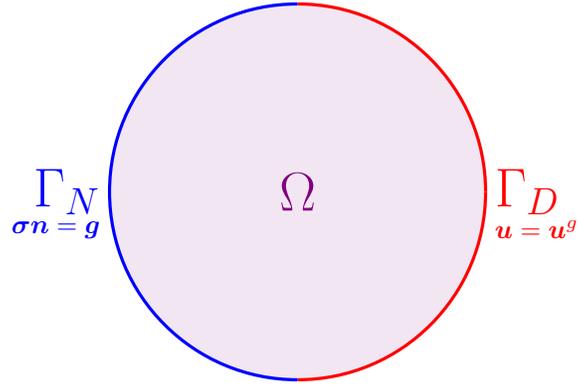
\begin{figure}[tbp]
\centering
\begin{tikzpicture}[scale = 2.5]
  \filldraw[color=white, fill=violet!10,  thick](0,0) circle (1);
  \draw [red, very thick] (0:1) arc [radius=1, start angle=0, end angle=90];
  \draw [red, very thick] (270:1) arc [radius=1, start angle=270, end angle=360];
  \draw [blue, very thick] (90:1) arc [radius=1, start angle=90, end angle=270];
\draw(0,0)[color = violet!100]node{\huge$\Omega$};
\draw(1,0)[right, color = red!100]node{\huge$\Gamma_D$};
\draw(-1,0)[left, color = blue!100]node{\huge$\Gamma_N$};
\draw(-1,-0.2)[left, color = blue!100]node{$\mbf{\sigma}\mbf{n} = \mbf{g}$};
\draw(1,-0.2)[right, color = red!100]node{$\mbf{u} = \mbf{u}^g$};
\end{tikzpicture}
\caption{Test case with mixed boundary conditions: the geometry of Dirichlet and Neumann boundary parts.}\label{fig:illustration_mixed_elasticity}
\end{figure}

\begin{figure}[tbp]
\centering 
 \includegraphics[width = 0.45\textwidth]{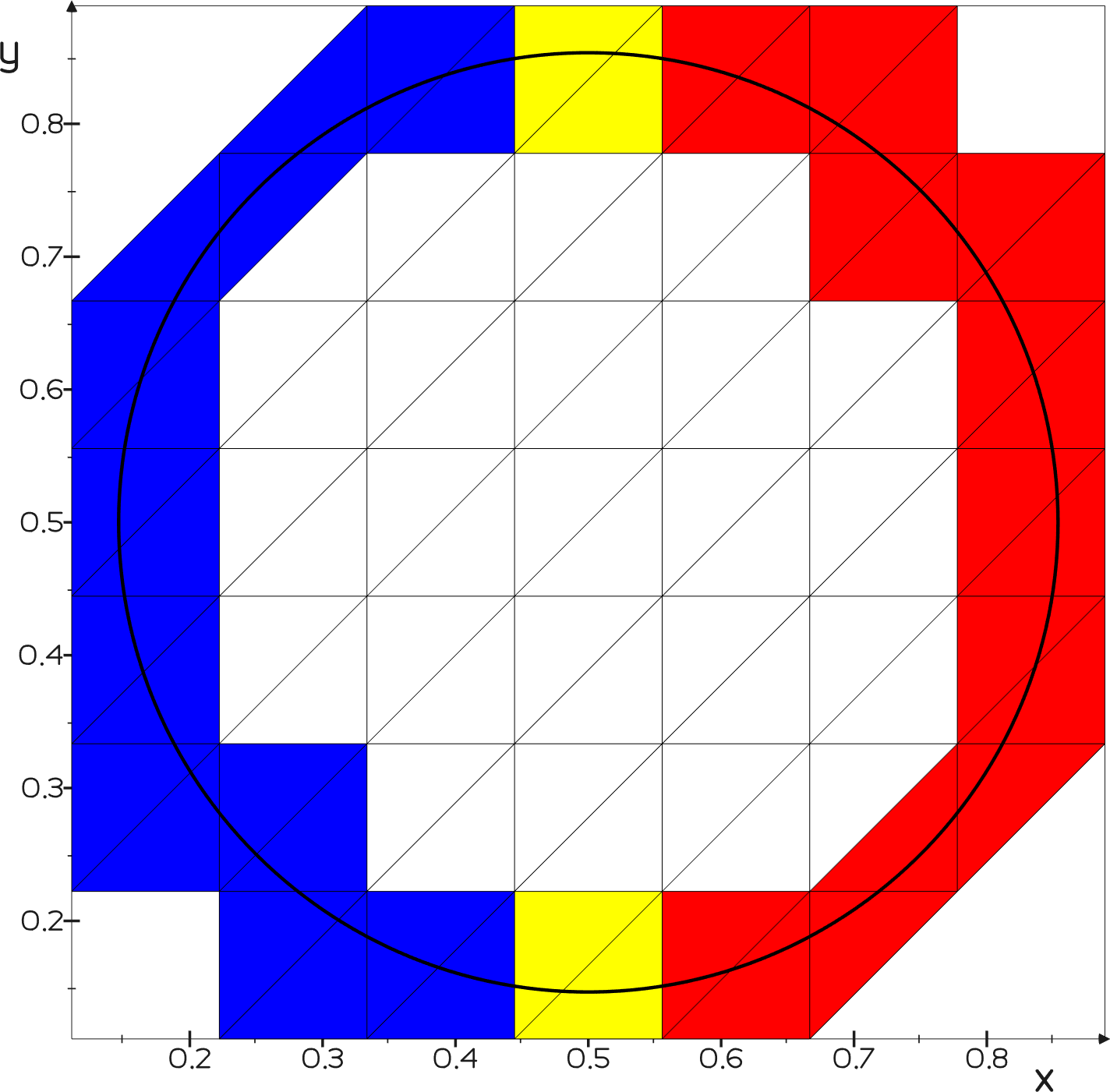}
\quad \includegraphics[width=0.45\textwidth]{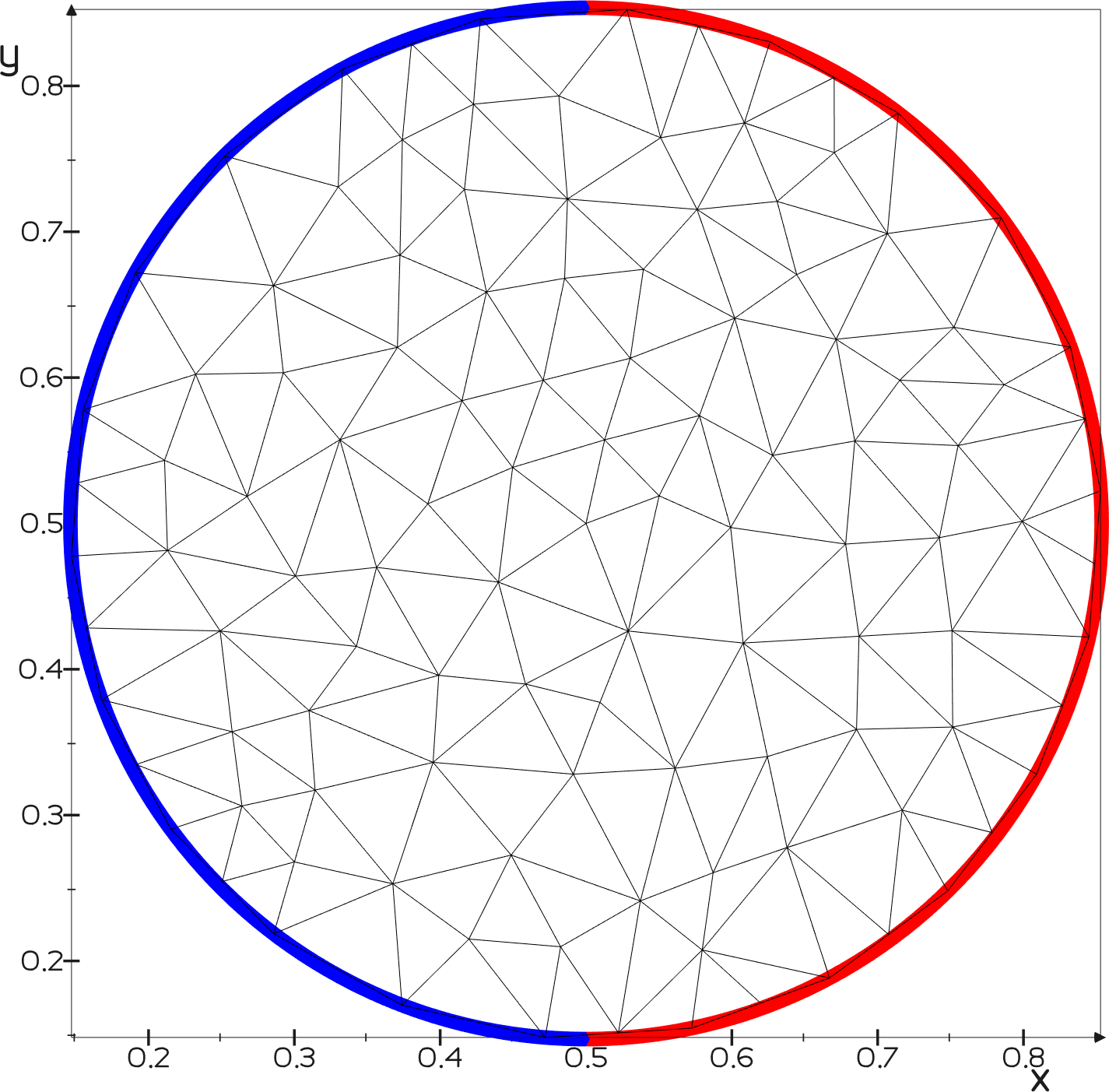}
\caption{Test case with mixed boundary conditions, meshes not resolving the Dirichlet/Neumann junction. Left: active meshes for $\phi$-FEM, red for $\Th^{\Gamma_D}$, blue  for $\Th^{\Gamma_N}$, yellow for $\Th^{\Gamma}$ otherwise unmarked. Right: a mesh for standard FEM, red boundary facets on ${\Gamma_D}$, blue boundary facets on ${\Gamma_N}$, note that some boundary facets contain both Dirichlet and Neumann parts. }\label{fig:meshes mixed 2}
\end{figure}

To set up the geometry of the problem, we recall that the domain $\Omega$ is given by the level set function $\phi$, $\Omega=\{\phi<0\}$, and assume furthermore that the boundary  partition into the Dirichlet and Neumann parts is governed by a secondary level set $\psi$, $$\Gamma_D=\Gamma\cap\{\psi<0\},\quad \Gamma_N=\Gamma\cap\{\psi>0\}\,.$$ Introducing the active meshes $\Th$ and $\Th^\Gamma$ as above, cf. (\ref{Th}), (\ref{ThGam}), and Fig.~\ref{fig}, we want now further partition the submesh $\Th^\Gamma$ into two parts: $\Th^{\Gamma_D}$ around $\Gamma_D$, serving to impose the Dirichlet boundary conditions, and $\Th^{\Gamma_N}$ around $\Gamma_N$ for the Neumann ones. The natural choice for these is
\begin{equation}\label{ThGamDN}
    \mathcal{T}_h^{\Gamma_D} := \{ T \in \mathcal{T}_h^{\Gamma} : \psi\leqslant 0\text{ on }T \} \qquad \text{ and } \qquad 
\mathcal{T}_h^{\Gamma_N} := \{ T \in \mathcal{T}_h^{\Gamma} : \psi\geqslant 0\text{ on }T \}\,.
\end{equation} 
As before, we denote the domains occupied by meshes $\Th$,$\mathcal{T}_h^{\Gamma}$,$\mathcal{T}_h^{\Gamma_D}$,$\mathcal{T}_h^{\Gamma_N}$ by $\Omega_h$,$\Omega_h^\Gamma$,$\Omega_h^{\Gamma_D}$,$\Omega_h^{\Gamma_N}$ respectively.
Note that these definitions may leave a small number of cells of $\mathcal{T}_h^{\Gamma}$ out of both $\mathcal{T}_h^{\Gamma_D}$ and $\mathcal{T}_h^{\Gamma_N}$. Indeed, there may be mesh cells, near the junction of Dirichlet and Neumann parts, where $\psi$ changes sign inside the cell, so that $\psi$ is neither everywhere positive not everywhere negative on such a cell.   This is illustrated at Fig.~\ref{fig:meshes mixed 2} (left) where the Dirichlet/Neumann junction is supposed at $x=0.5$, i.e. the secondary level set is $\psi(x,y)=0.5-y$, c.f. Fig.~\ref{fig:illustration_mixed_elasticity}. The active mesh cells intersected by $\Gamma$ on Fig.~\ref{fig:meshes mixed 2} are either on the Dirichlet side (they form thus $\Th^{\Gamma_D}$  and are colored in red), or on on the Neumann side (they form thus $\Th^{\Gamma_N}$  and are colored in blue), or in between (they are then in $\Th^{\Gamma}$ but not in $\Th^{\Gamma_D}$ or $\Th^{\Gamma_N}$, and are colored in yellow).

Assuming once more that $\mbf{u}$, the solution to (\ref{eq:elast})--(\ref{bcGammaD})--(\ref{bcGammaN}), can be extended from $\Omega$ to $\Omega_h$ as the solution to  the same governing equation (\ref{eq:elast}), we introduce a $\phi$-FEM scheme, combining the Dual $\phi$-FEM Dirichlet approach, as introduced in (\ref{DualPhi}) and (\ref{DualDiscrete}), with the indirect imposition of Neumann boundary condition as proposed in \cite{phiFEM2}. We thus keep $\mbf{u}$ as the primary unknown on $\Omega_h$ and recall that it satisfies the variational formulation (\ref{weakExt}). The Dirichlet boundary condition affects the solution on  $\Omega_h^{\Gamma_D}$ through the introduction of the auxiliary variable $\mbf{p}_D$ there. We thus adapt (\ref{ansatzDual}) from the pure Dirichlet case as  \begin{equation}\label{ansatzDirMixt}
  	\mbf{u}=\mbf{u}^g + \phi \mbf{p}_D,\text{ on }\Omega_h^{\Gamma_D}\,.
\end{equation} 
We have assumed here that $\mbf{u}^g$ is extended from ${\Gamma_D}$ to $\Omega_h^{\Gamma_D}$.

The Neumann boundary condition will affect $\mbf{u}$ on  $\Omega_h^{\Gamma_N}$ through the introduction of two auxiliary variables there. We first introduce a tensor-valued variable $\mbf{y}$ on $\Omega_h^{\Gamma_N}$ setting $\mathbf{y}=-\mbf{\sigma}(\mbf{u})$. It remains to impose $\mbf{y}\mbf{n}=-\mbf{g}$ on $\Gamma_N$. To this end, we note that the outward-looking unit normal $\mbf{n}$ is given on $\Gamma$ by $\mbf{n}=\frac{1}{|\nabla\phi|}\nabla\phi$ so that the Neumann boundary condition is satisfied by setting $\mbf{y}\nabla\phi+\mbf{g}|\nabla\phi| =- \mbf{p}_N\phi$ on $\Omega_h^{\Gamma_N}$ where $\mbf{p}_N$ is yet another (vector-valued) auxiliary variable on $\Omega_h^{\Gamma_N}$. This can be summarized as
	\begin{subequations}\label{phiN}
	\begin{align}
	\mathbf{y}+\mbf{\sigma}(\mbf{u})=0,&\quad \text{on } \Omega_h^{\Gamma_N}\,, \\
	\mbf{y}\nabla\phi+ \mbf{p}\phi=-\mbf{g}|\nabla\phi|,&\quad \text{on } \Omega_h^{\Gamma_N}\,.
	\end{align}
	\end{subequations}
Note that the combination of (\ref{ansatzDirMixt}) with (\ref{phiN}a-b) does not impose the mixed Dirichlet/Neumann conditions on the whole of $\Gamma$ since the latter may be not completely covered by  $\Omega_h^{\Gamma_D}\cup\Omega_h^{\Gamma_N}$. Fortunately, this defect of the formulation on the continuous level can be repaired on the discrete level by adding the appropriate stabilization to the FE discretization.  	

To describe the resulting FE scheme, we start by introducing the FE spaces. As before, we fix an integer $k\ge 1$ and keep the space $V_h$, as  defined in (\ref{spaceVh}), for the approximation $\mbf{u}_h$ of the primary variable $\mbf{u}$. We need also the spaces for the  approximation of the auxiliary variables $\mbf{p}_{h,D}$ and  $\mbf{p}_{h,N}$, respectively $Q_h^k(\Omega_h^{\Gamma_D})$ and $Q_h^{k-1}(\Omega_h^{\Gamma_N})$ as defined in (\ref{spaceQh}), as well as the space $Z_h(\Omega_h^{\Gamma_N})$ to approximate $\mbf{y}$,
where for each submesh $\mathcal{M}_h$ of $\Th$, $Z_h(\mathcal{M}_h)$ is defined by
\begin{equation}
\label{espaceZh}
 Z_{h}(\mathcal{M}_h) := \left\lbrace \mbf{z}_h :\mathcal{M}_h \to\mathbb{R}^{(d\times d) } : \mbf{z}_{h |T} \in \mathbb{P}^k(T)^{(d\times d)} \ \ \forall T \in \mathcal{M}_h, \ \mbf{z}_h \text{ continuous on }\mathcal{M}_h \right\rbrace \,.
 \end{equation}  
Now, combining the variational formulation (\ref{weakExt}) with the (\ref{ansatzDirMixt}) and (\ref{phiN}a-b) imposed in a least-squares manner, we get the following scheme: find $\mbf{u}_h \in V_h$,  $\mbf{p}_{h,D} \in Q_{h}^k(\Omega_h^{\Gamma_D})$, $\mbf{y}_h \in Z_h(\Omega_h^{\Gamma_N})$ and $\mbf{p}_{h,N} \in Q_{h}^{k-1}(\Omega_h^{\Gamma_N})$ such that  
\begin{multline}\label{mixtDicrete}
\int_{\Omega_h} {\mbf{\sigma}}(\mbf{u}_h) : \nabla\mbf{v}_h 
 - \int_{\partial\Omega_h\setminus\partial\Omega_{h,N}}\mbf{\sigma}(\mbf{u}_h)\mbf{n} \cdot \mbf{v}_h  
 + \int_{\partial\Omega_{h,N}}\mbf{y}_h\mbf{n} \cdot \mbf{v}_h 
 \\ + \gamma_u \int_{\Omega_h^{\Gamma_N}}(\mbf{y}_h+ {\mbf{\sigma}}(\mbf{u}_h)) : (\mbf{z}_h+ {\mbf{\sigma}}(\mbf{v}_h))   + \frac{\gamma_p}{h^2} \int_{\Omega_h^{\Gamma_N}}\left( \mbf{y}_h\nabla \phi_h + \frac{1}{h}\mbf{p}_{h,N} \phi_h \right)\cdot\left( \mbf{z}_h\nabla \phi_h + \frac{1}{h}\mbf{q}_{h,N}\phi_h \right)
 \\ + \frac{\gamma}{h^2} \int_{\Omega_h^{\Gamma_D}} (\mbf{u}_h - \frac{1}{h}\phi_h \mbf{p}_{h,D}) \cdot (\mbf{v}_h - \frac{1}{h}\phi_h \mbf{q}_{h,D}) 
 + G_h(\mbf{u}_h,\mbf{v}_h) 
 + J_h^{lhs,D}(\mbf{u}_h,\mbf{v}_h) + J_h^{lhs,N}(\mbf{y}_h,\mbf{z}_h) \\ 
= \int_{\Omega_h} \mbf{f} \cdot \mbf{v}_h  + \frac{\gamma}{h^2} \int_{\Omega_h^D} \mbf{u}^g_h \cdot (\mbf{v}_h - \frac{1}{h}\phi_h \mbf{q}_{h,D})  - \frac{\gamma_p}{h^2} \int_{\Omega_h^{\Gamma_N}} \mbf{g} \cdot |\nabla \phi_h| (\mbf{z}_h \cdot \nabla \phi_h + \frac{1}{h} \mbf{q}_{h,N} \phi_h)  \\
+ J_h^{rhs,D}(\mbf{v}_h)  + J_h^{rhs,N}(\mbf{z}_h), \\ \forall \mbf{v}_h\in V_h, \mbf{q}_{h,D} \in Q_{h}^k(\Omega_h^{\Gamma_D}), \mbf{z}_h \in Z_h(\Omega_h^{\Gamma_N}), \mbf{q}_{h,N} \in Q_{h}^{k-1}(\Omega_h^{\Gamma_N})\,.
\end{multline}
We have added here the ghost stabilization $G_h$ defined by (\ref{Gh}) as in the pure Dirichlet case.
The additional stabilizations terms $J_h^{lhs}$, $J_h^{rhs}$ are now adapted from (\ref{Jhlrhs}) and separated into the terms acting on $\mbf{u}_h$ on  the Dirichlet cells of $\Th^\Gamma$ (and also those not marked), and the terms acting on $\mbf{y}_h$ on  the Neumann cells:   
\begin{align}
J_h^{lhs,D}(\mbf{u}, \mbf{v}) &: =  \sigma_D h^2  \sum_{T \in \Th^{\Gamma} \setminus \Th^{\Gamma_N}} \int_T  \Div \mbf{\sigma}(\mbf{u})\cdot \Div\mbf{\sigma}(\mbf{v}) \,,
&& \notag \\
J_h^{rhs,D} (\mbf{v}) &: =- \sigma_D h^2\sum_{T \in \Th^{\Gamma} \setminus \Th^{\Gamma_N}} \int_T \mbf{f} \cdot\Div\mbf{\sigma}(\mbf{v}) 
\,,
&& \notag \\
    J_h^{lhs,N}(\mbf{y},\mbf{z}) &=\gamma_{div} \int_{\Omega_h^{\Gamma_N}}\Div \mbf{y} \cdot \Div \mbf{z}
\,,
& J_h^{rhs,N}(\mbf{z}) &= \gamma_{div} \int_{\Omega_h^{\Gamma_N}} \mbf{f} \cdot \Div\mbf{z}.
\label{JhN}
\end{align}
These stabilizations are consistent with the governing equations $\Div \mbf{\sigma}(\mbf{u})=-\mbf{f}$, rewritten as  $\Div \mbf{y}=\mbf{f}$, using (\ref{phiN}a), wherever possible, i.e. on $\Omega_h^{\Gamma_N}$. Note that a similar treatment is applied to the boundary integral terms on $\partial\Omega$ in (\ref{weakExt}). In \eqref{mixtDicrete}, they are rewritten  in terms of $\mbf{y}$, using (\ref{phiN}a)  and (\ref{phiN}b), wherever possible. We thus introduce a part of the boundary $\partial\Omega_h$, referred to as $\partial\Omega_{h,N}$, formed by the boundary facets of $\Th$  belonging to the cells in $\Th^{\Gamma_N}$. We replace ${\mbf{\sigma}}(\mbf{u}_h)$ by $-\mbf{y}_h$ on $\partial\Omega_{h,N}$, while keeping the boundary term as is on the remaining part of the boundary. All this contributes to the coerciveness of the bilinear form in (\ref{mixtDicrete})  and good conditioning of the matrix as can be proven following the ideas of \cite{phiFEM2}.  We emphasize again that neither Dirichlet nor Neumann boundary conditions are imposed in any way in scheme (\ref{mixtDicrete}) on the cells in $\Th^{\Gamma}\setminus(\Th^{\Gamma_D}\cup\Th^{\Gamma_N})$ (the cells in yellow on Fig.~\ref{fig:meshes mixed 2}). On the other hand, both stabilizations $G_h$ and $J_h$ are active on the whole $\Th^{\Gamma}$, comprising these cells not marked as Dirichlet or Neumann.

\begin{figure}[p]
\centering 
 \includegraphics[width = 0.45\textwidth]{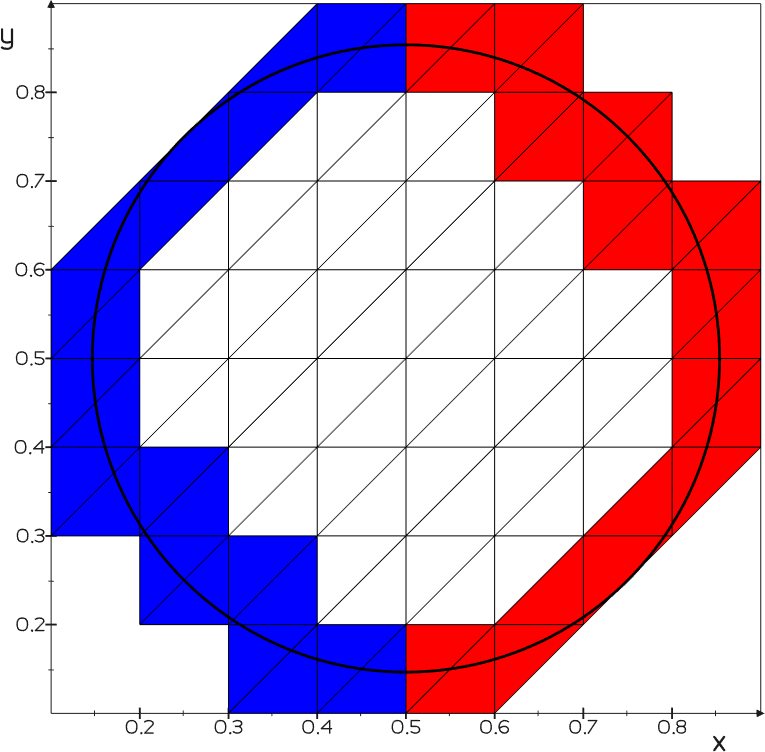}
\quad \includegraphics[width=0.45\textwidth]{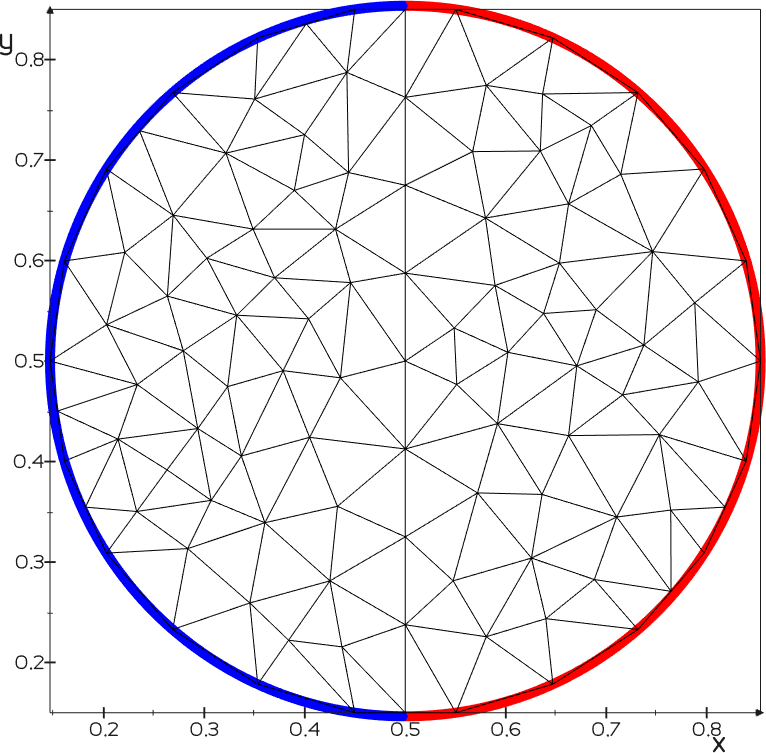}
\caption{Test case with mixed boundary conditions, meshes resolving the Dirichlet/Neumann junction. Left: active meshes for $\phi$-FEM, red for $\Th^{\Gamma_D}$, blue  for $\Th^{\Gamma_N}$. Right: a mesh for standard FEM, red boundary facets on ${\Gamma_D}$, blue boundary facets on ${\Gamma_N}$. }\label{fig:meshes mixed 1}
\end{figure}

\begin{figure}[p]
\centering
\begin{tikzpicture}
\begin{loglogaxis}[name = ax1, width = .45\textwidth, xlabel = $h$, 
            ylabel = relative error,
            legend style = { at={(0.8,1.35)}, legend columns =1,
			/tikz/column 2/.style={column sep = 10pt}}]
\addplot coordinates {
(0.08838834764831845,0.00019162789405352987)
(0.04419417382415922,9.019019807546487e-06)
(0.02209708691207961,7.722996330926221e-07)
(0.011048543456039806,6.463963785481317e-08)
(0.005524271728019903,6.723748171411495e-09)};
\addplot coordinates { 
(0.08838834764831845,0.0022236508129968484)
(0.04419417382415922,0.00025661468575222617)
(0.02209708691207961,4.506471015826728e-05)
(0.011048543456039806,8.087493674563759e-06)
(0.005524271728019903,1.6960635243482559e-06)};
\addplot coordinates { 
(0.08789154903752427,0.0035246050135470824)
(0.044129789705730664,0.0008236790908525067)
(0.022096512804264105,0.0002031809648554605)
(0.011048263536162476,5.020191340038304e-05)
(0.005524239072190211,1.2340916877939585e-05)};
\addplot coordinates{
(0.08789154903752427,0.016974433240575463)
(0.044129789705730664,0.006034873183160213)
(0.022096512804264105,0.0015285621307403931)
(0.011048263536162476,0.0007285213745054552)
(0.005524239072190211,0.0002601134487378277)};
\logLogSlopeTriangle{0.53}{0.2}{0.45}{2}{red};
\logLogSlopeTriangle{0.53}{0.2}{0.12}{3}{blue};
\legend{$L^2$ error $\phi$-FEM, $H^1$ error $\phi$-FEM, $L^2$ error standard FEM, $H^1$ error standard FEM}
\end{loglogaxis}
\begin{loglogaxis}[width = .45\textwidth, ylabel = Computing time (s), name = ax2, at = {(ax1.south east)}, xshift = 2cm,
            xlabel =$L^2$ relative error, 
            legend style = { at={(1,1.2)}, legend columns =2,
			/tikz/column 2/.style={column sep = 10pt}}]
\addplot coordinates {
(0.00019162789405352987,0.07287120819091797)
(9.019019807546487e-06,0.13202333450317383)
(7.722996330926221e-07,0.3721940517425537)
(6.463963785481317e-08,1.2812981605529785)
(6.723748171411495e-09,5.520588159561157)};
\addplot coordinates { 
(0.0035246050135470824,0.020807504653930664)
(0.0008236790908525067,0.05296206474304199)
(0.0002031809648554605,0.21299290657043457)
(5.020191340038304e-05,0.9861962795257568)
(1.2340916877939585e-05,4.93593955039978)};
\legend{$\phi$-FEM, Standard FEM}
\end{loglogaxis}
\end{tikzpicture}
\caption{Test case with mixed boundary conditions, results on meshes as on Fig.~\ref{fig:meshes mixed 1}. Left: $L^2$ and $H^1$ relative errors under the mesh refinement.  Right: computing time vs. the $L^2$ relative error. }\label{fig:mixed_elasticity}
\end{figure}
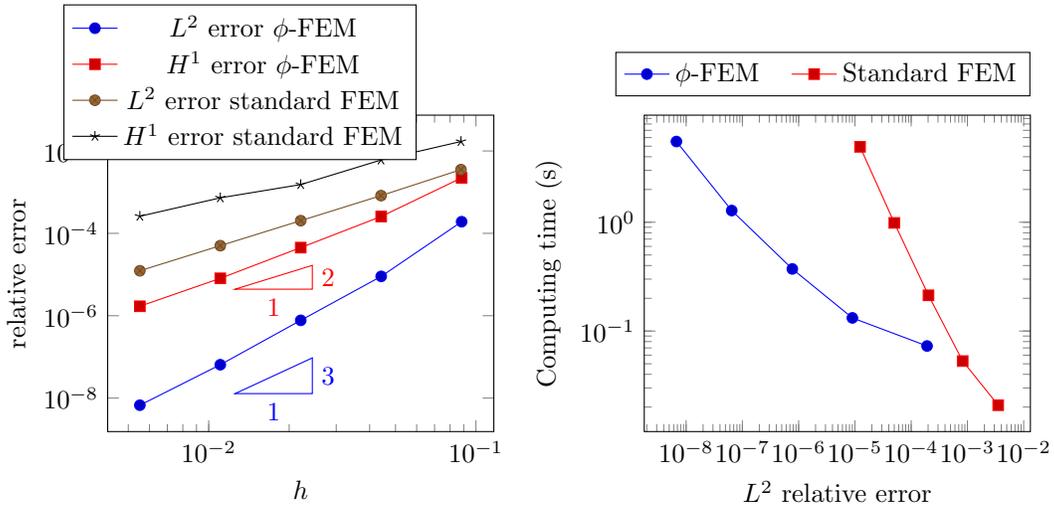

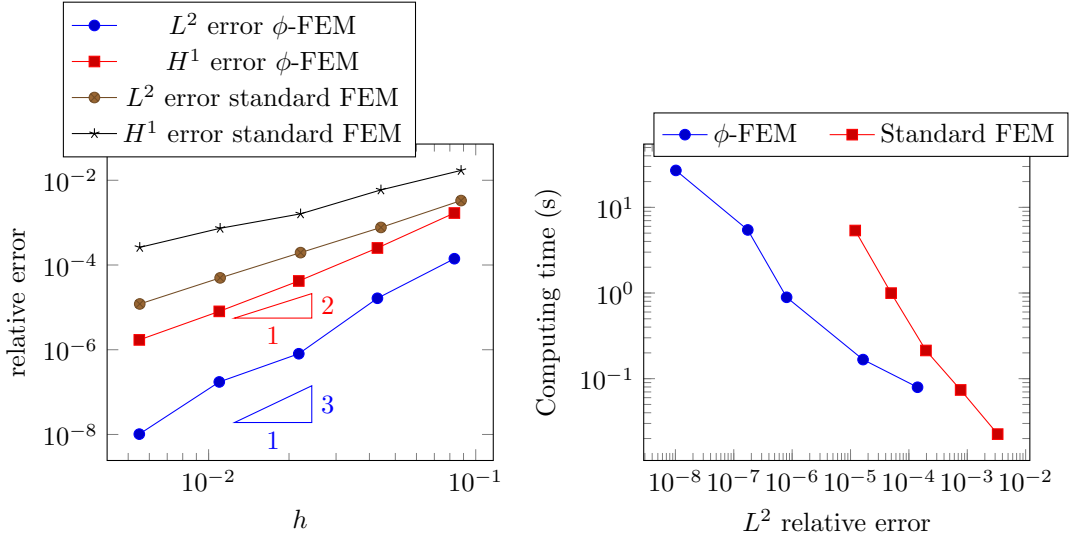
\begin{figure}[tbp]
\centering
\begin{tikzpicture}
\begin{loglogaxis}[name = ax1, width = .45\textwidth, xlabel = $h$, 
ylabel = relative error,
legend style = { at={(0.8,1.45)}, legend columns =1,
			/tikz/column 2/.style={column sep = 10pt}}]
\addplot coordinates {
(0.08318903308077032,0.00014071942503931245)
(0.04285495643554845,1.6398220861748423e-05)
(0.021757131728816926,8.037137350273614e-07)
(0.010962895832349594,1.7379414354329388e-07)
(0.005502776507288445,1.0190675532414034e-08)};
\addplot coordinates { 
(0.08318903308077032,0.001685377218347434)
(0.04285495643554845,0.00025228072880838193)
(0.021757131728816926,4.1981937129849344e-05)
(0.010962895832349594,8.04580293961358e-06)
(0.005502776507288445,1.7062245894835732e-06)};
\addplot coordinates { 
(0.08814237266741118,0.003318016542529396)
(0.044176082108213596,0.0007701415118214259)
(0.022095901417288732,0.0001959052803742703)
(0.01104624219043868,4.965445115252627e-05)
(0.005523930001854928,1.2033909691066592e-05)};
\addplot coordinates{
(0.08814237266741118,0.0170779277806841)
(0.044176082108213596,0.005916646533061518)
(0.022095901417288732,0.001611950067370548)
(0.01104624219043868,0.0007352560307018717)
(0.005523930001854928,0.000260614890322195)};
\logLogSlopeTriangle{0.53}{0.2}{0.45}{2}{red};
\logLogSlopeTriangle{0.53}{0.2}{0.12}{3}{blue};
\legend{$L^2$ error $\phi$-FEM, $H^1$ error $\phi$-FEM, $L^2$ error standard FEM, $H^1$ error standard FEM}
\end{loglogaxis}
\begin{loglogaxis}[width = .45\textwidth, ylabel = Computing time (s), name = ax2, at = {(ax1.south east)}, xshift = 2cm,
xlabel =$L^2$ relative error, 
legend style = { at={(1.1,1.1)}, legend columns =2,
/tikz/column 2/.style={column sep = 10pt}}]
\addplot coordinates {
(0.00014071942503931245,0.07939004898071289)
(1.6398220861748423e-05,0.16697168350219727)
(8.037137350273614e-07,0.8911266326904297)
(1.7379414354329388e-07,5.4493019580841064)
(1.0190675532414034e-08,26.988556385040283)};
\addplot coordinates { 
(0.003318016542529396,0.022522687911987305)
(0.0007701415118214259,0.07380819320678711)
(0.0001959052803742703,0.21341729164123535)
(4.965445115252627e-05,0.9971983432769775)
(1.2033909691066592e-05,5.366613149642944)};
\legend{$\phi$-FEM, Standard FEM}
\end{loglogaxis}
\end{tikzpicture}
\caption{Test case with mixed boundary conditions, results on meshes as on Fig.~\ref{fig:meshes mixed 2}. Left: $L^2$ and $H^1$ relative errors under the mesh refinement.  Right: computing time vs. the $L^2$ relative error. }\label{fig:mixed_elasticity_2}
\end{figure}

\paragraph*{Test case: }
We are now going to present some numerical results with method (\ref{mixtDicrete}) highlighting the optimal convergence of $\phi$-FEM and comparing it with a standard FEM.
We use the same geometry (\ref{phiCircle}), elasticity parameters and the exact solution (\ref{manufu}) as for the case of pure Dirichlet conditions on page~\pageref{testcaseDir}. We set furthermore the Dirichlet boundary conditions (\ref{bcGammaD}) for $x>0.5$ and the Neumann boundary conditions (\ref{bcGammaN}) for $x<0.5$, c.f. Fig.~\ref{fig:illustration_mixed_elasticity}, i.e. we choose the secondary level set as $\psi=0.5-x$. The data $\mbf{u}^g$ and $\mbf{g}$ are computed from the exact solution. In $\phi$-FEM they should be extended from $\Gamma$ to appropriate portion of the strip $\Omega_h^\Gamma$. We choose these extensions as
\[\begin{cases} \mbf{u}^g = \mbf{u}_{ex}(1+\phi), \quad &\text{ on } \Omega_h^{\Gamma} \cap \{ x \geqslant 0.5 \} \,, \\ 
\mbf{g} = \mbf{\sigma}(\mbf{u}_{ex}) \frac{\nabla \phi}{\| \nabla \phi\|} + \mbf{u}_{ex}  \phi, \quad &\text{ on }  \Omega_h^{\Gamma} \cap \{ x < 0.5 \}\,. \end{cases}  \]
Again, both expressions are perturbed away from  $\Gamma$ to mimic the real-life situation where the data are available only on $\Gamma$. The stabilization parameters are set to $\gamma_{div} = \gamma_u = \gamma_p = 1.0$, $\sigma = 0.01$  and $\gamma = \sigma_D = 20.0$. 

We start by studying mesh configurations where the Dirichlet-Neumann junction line $\{x=0.5\}$ happens to be covered by the mesh facets both in the background mesh used by $\phi$-FEM, and in the fitted mesh used by FEM, as illustrated in Fig.~\ref{fig:meshes mixed 1}. All the boundary cells in $\Th^\Gamma$ are marked in this case either as Dirichlet or as Neumann ones, according to the criterion (\ref{ThGamDN}), giving, respectively, red and blue cells on Fig.~\ref{fig:meshes mixed 1} (left). There is no ambiguity for the standard FEM fitted meshes: all the boundary facets are straightforwardly marked either as Dirichlet or as Neumann, cf. Fig.~\ref{fig:meshes mixed 1} (right) with the same color code as for the unfitted mesh. The results obtained by both $\phi$-FEM (\ref{mixtDicrete}) and the standard FEM, using $\mathbb{P}^2$-Lagrange polynomials for $\mbf{u}_h$ in both cases, are reported in Fig.~\ref{fig:mixed_elasticity}. On the left, the relative errors are plotted with respect to the mesh step. We observe again the optimal convergence orders for $\phi$-FEM, while the convergence of the standard FEM is sub-optimal in the $L^2$-norm. The $\phi$-FEM approach is again systematically more precise in both norms. On the right side of the same figure, we plot the computing times and notice again that $\phi$-FEM is less expensive than the standard FEM.


Let us now turn to a less artificial mesh configuration where the Dirichlet/Neumann junction point can turn up inside a mesh cell of the background mesh, or inside a boundary facet of the fitted mesh. We study these situations on a series of meshes, as illustrated in Fig.~\ref{fig:meshes mixed 2}. In the case of the background meshes used for $\phi$-FEM, we ensure in particular that there are no vertical grid line with the abscissa $x=0.5$ so that there are exactly 4 cells cells in $\mathcal{T}_h^{\Gamma}$ that are neither  in $\mathcal{T}_h^{\Gamma_D}$ nor in $\mathcal{T}_h^{\Gamma_N}$ (yellow cells on the left side of Fig.~\ref{fig:meshes mixed 2}). We recall that scheme (\ref{mixtDicrete}) does not impose any boundary conditions on these cells, but retains the stabilization there (in particular, the governing equation is still re-enforced on these cells in the least squares manner). Note that the fitted FEM is not straightforward to implement in this case either, since the Dirichlet boundary conditions cannot be strongly imposed  on the boundary facets which turn up only partially on the Dirichlet side. We bypass this difficulty by treating the Dirichlet conditions by penalization, so that the "standard" FEM is now defined as: find $\mbf{u}_h$ in the $\mathbb{P}^k$ FE space (without any restrictions on the boundary) such that 
\begin{equation}\label{FEMmixed}
  \int_{\Omega} \mbf{\sigma} (\mbf{u}_h) : \nabla\mbf{v}_h   
    + \frac{1}{\varepsilon}\int_{\Gamma_D} \mbf{u}_h \cdot \mbf{v}_h
   = \int_{\Omega} \mbf{f} \cdot \mbf{v}_h+
   \int_{\Gamma_N} \mbf{g} \cdot \mbf{v}_h
    + \frac{1}{\varepsilon}\int_{\Gamma_D} \mbf{u}^g \cdot \mbf{v}_h
\end{equation}
for all $\mbf{v}_h$ in the same FE space as $\mbf{u}_h$, with a small parameter $\varepsilon>0$.

The mesh refinement study in this case is reported at Fig.~\ref{fig:mixed_elasticity_2}. Comparing the results with those of Fig.~\ref{fig:mixed_elasticity} (obtained on idealized unrealistic meshes without any unmarked cells), we observe that the behavior of $\phi$-FEM (\ref{mixtDicrete}) is almost unaffected by the presence (or not) of the unmarked ``yellow" cells, although the convergence curve for the $L^2$ relative error is now slightly less regular. In particular, the conclusions about the relative merits of $\phi$-FEM and the fitted  FEM, now in version (\ref{FEMmixed}),  remain unchanged: $\phi$-FEM is  more precise on comparable meshes and less expensive in terms of the computing times for a given error tolerance. 

\FloatBarrier

\section{Linear elasticity with multiple materials.}\label{sectInterface}

We now consider the case of interfaces problems, i.e. partial differential equations  with coefficients jumping across an interface, which can cut the computational mesh in an arbitrary manner. The simplest meaningful example in the realm of linear elasticity is given by structures consisting of multiple materials having different elasticity parameters. 
This situation has already been treated in XFEM \cite{carraroXFEMinterface, NitscheXFEMinterface, xiaoXFEM, Xiao}, CutFEM \cite{CutFEMinterface, Hansbo2005interface,hansbointerface,lehrenfeldinterface}, and SBM \cite{interfaceSBMli} paradigms. 
We are now going to demonstrate the applicability of $\phi$-FEM in this context.

Let us assume that the structure occupies a domain $\Omega$ and it consists of two materials that occupy two subdomains $\Omega_1$ and $\Omega_2$ separated by the interface $\Gamma$.  To fix the ideas, we further assume that $\Omega_1$ is surrounded by $\Omega_2$, so that the interface $\Gamma$ can actually be described as  $\Gamma = \partial\Omega_1$, as illustrated at Fig.~\ref{fig:interface}. We also assume that the displacement $\mbf{u}$ is given on the external boundary (these assumptions are not restrictive and the forthcoming method can be easily adapted to other situations, e.g. with $\Gamma$ touching $\partial\Omega$ or with Neumann boundary conditions on the external boundary). We then consider the problem for the displacement $\mbf{u}$ on $\Omega$: 
\begin{equation}\label{eq:interface_elasticity}
\begin{cases}
- \Div \mbf{\sigma} (\mbf{u}) &= \mbf{f} \,, \text{ on } \ \Omega\backslash\Gamma \,, \\
\mbf{u} &= \mbf{u}^g \,, \text{ on } \ \partial \Omega \,, \\
[\mbf{u} ] &= 0\,, \text{ on } \ \Gamma \,, \\
[\mbf{\sigma}(\mbf{u})  \mbf{n}] &= 0 \,, \text{ on } \ \Gamma \,, 
\end{cases}
\end{equation}
where $\mbf{n}$ is the unit normal pointing from $\Omega_1$ to $\Omega_2$, and the brackets $[\cdot]$ stand for the jump across $\Gamma$. 
The elasticity parameters are assumed constant on each sub-domain, but different from each other. The stress tensor is thus given by
\[ \mbf{\sigma} (\mbf{u}) =
\begin{cases}
\mbf{\sigma}_1(\mbf{u})=2\mu_1 \mbf{\varepsilon}(\mbf{u})+\lambda_1 (\Div \mbf{u})I\,,\text{ on }\Omega_1\,,  \\ 
\mbf{\sigma}_2(\mbf{u})=2\mu_2 \mbf{\varepsilon}(\mbf{u})+\lambda_2 (\Div \mbf{u})I\,,\text{ on }\Omega_2\,,
\end{cases}
\]
with  the Lamé parameters $\lambda_i $ and $\mu_i$ defined via the formulas (\ref{muEnu}) with given $E_i,\nu_i$, $i=1,2$. Introducing the displacements $\mbf{u}_i=\mbf{u} |_{\Omega_i}$, $i=1,2$ on $\Omega_1$ and $\Omega_2$ separately, problem (\ref{eq:interface_elasticity}) can be rewritten as the system of two coupled sub-problems: 
\begin{equation}\label{eq:interface_elasticity2}
\begin{cases}
- \Div \mbf{\sigma}_i (\mbf{u}_i) &= \mbf{f} \,, \text{ on } \ \Omega_i\,,\ i=1,2, \\
\mbf{u}_2 &= \mbf{u}^g \,, \text{ on } \ \partial \Omega \,, \\
\mbf{u}_1  &= \mbf{u}_2\,, \text{ on } \ \Gamma \,, \\
\mbf{\sigma}_1(\mbf{u}_1) \mbf{n} &= \mbf{\sigma}_2(\mbf{u}_2) \mbf{n} \,, \text{ on } \ \Gamma \,.
\end{cases}
\end{equation}
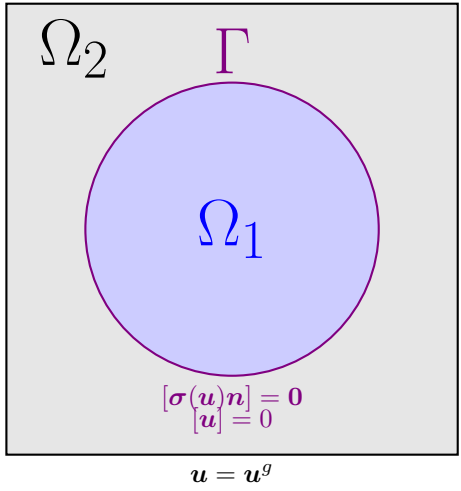
\begin{figure}[htbp]
\centering
\begin{tikzpicture}[scale=3]
\filldraw[color=black!100, fill=lightgray!40, thick]rectangle (2,2);
\filldraw[color=violet!100, fill=blue!20,  thick](1,1) circle (0.65);
\draw(1,1)[color = blue!100]node{\Huge$\Omega_1$};
\draw(1,1.65)[above, color = violet!100]node{\Huge$\Gamma$};
\draw(1,0.35)[below, color = violet!100]node{$[\mbf{\sigma}(\mbf{u}) \mbf{n} ] = \mbf{0}$};
\draw(1,0.25)[below, color = violet!100]node{$[\mbf{u}] = 0$};
\draw(0.3,1.80)[color = black!100]node{\Huge$\Omega_2$};
\draw(1,0.0)[below, color = black!100]node{$\mbf{u} = \mbf{u}^g$};
\end{tikzpicture}
\caption{Geometry with the interface $\Gamma$: elasticity with multiple materials.}\label{fig:interface}
\end{figure}

We suppose that $\Omega$ is sufficiently simple-shaped so that a matching mesh $T_h$ on $\Omega$ is easily available (again, this assumption is not restrictive; we have seen that a complex-shape domain $\Omega$ can be also treated by $\phi$-FEM). On the contrary, the mesh $\Th$ is not supposed to match the internal interface $\Gamma$ and we are going to adapt $\phi$-FEM to this situation. The starting point is the reformulation (\ref{eq:interface_elasticity2}). We are thus going to discretize separately $\mbf{u}_1$ on $\Omega_1$ and $\mbf{u}_2$ on $\Omega_2$.  To this end, we introduce two active meshes $\mathcal{T}_{h,1}$ and $\mathcal{T}_{h,2}$, sub-meshes of $\Th$, constructed by retaining in $\mathcal{T}_{h,i}$ the cells of $\Th$ having a non-empty intersection with $\Omega_i$. In practice, the sub-domains are defined through a level-set $\phi$:
\begin{equation*}
    \Omega_1 = \{ \phi > 0 \}\cap\Omega, \qquad \Omega_2 = \{ \phi < 0 \}, \qquad \Gamma = \{ \phi = 0 \}\cap\Omega\, . 
\end{equation*} 
The sub-meshes $\mathcal{T}_{h,i}$ are defined using  a piecewise-polynomial approximation $\phi_h$ of $\phi$, rather than $\phi$ itself:
\begin{equation}\label{Th12}
\mathcal{T}_{h,1}:=\{T\in \Th:T\cap \{\phi_h>0\}\neq  \varnothing\}\text{ and }
\mathcal{T}_{h,2}:=\{T\in \Th:T\cap \{\phi_h<0\}\neq  \varnothing\}\,.
\end{equation}
We also introduce the sub-mesh $\Th^\Gamma$ as the intersection $\mathcal{T}_{h,1}\cap\mathcal{T}_{h,2}$ and denote by $\Omega_{h,1}$, $\Omega_{h,2}$, $\Omega_h^{\Gamma}$ the domains covered by meshes $\mathcal{T}_{h,1}$, $\mathcal{T}_{h,2}$, $\Th^\Gamma$ respectively. Similarly to  the simpler settings considered above, the unknowns  $\mbf{u}_1$ and $\mbf{u}_2$, living physically on $\Omega_1$ and $\Omega_2$,  will be discretized on larger domains $\Omega_{h,1}$ and  $\Omega_{h,2}$, introducing artificial extensions on narrow fictitious strips near $\Gamma$.  On the discrete level, the unknowns will be thus redoubled on the joint sub-mesh $\Th^\Gamma$. Several auxiliary unknowns will be introduced on $\Omega_h^\Gamma$ similar to the case of mixed boundary conditions above (indeed, we have to discretize both Dirichlet and Neumann conditions on the interface $\Gamma$ in the current setting). 

We now put the program above into the equations, first on the continuous level. Similarly to (\ref{weakExt}), the unknowns $\mbf{u}_i$ extended to larger domains $\Omega_h^i$ satisfy formally the variational formulations, cf. the first equation in (\ref{eq:interface_elasticity2}): 
\begin{equation}\label{weakExtomi} 
\int_{\Omega_{h,i}} \mbf{\sigma}_i (\mbf{u}_i) : \nabla\mbf{v}_i  
- \int_{\partial\Omega_{h,i}} \mbf{\sigma}_i (\mbf{u}_i)  \mbf{n}_i \cdot \mbf{v}_i
= \int_{\Omega_{h,i}} \mbf{f} \cdot \mbf{v}_i, \quad \forall \mbf{v}_i \text{ on }\Omega_{h_i}\text{ s.t.} \mbf{v}_i=\mbf{0}\text{ on }\partial\Omega\,.
\end{equation}
Here, with a slight abuse of notations, ${\partial\Omega_{h,i}}$ denotes the component of the boundary of $\Omega_{h,i}$ other than $\partial\Omega$, and  $\mbf{n}_i$ denotes the unit normal vector on $\partial\Omega_{h,i}$ pointing outside $\Omega_{h,i}$. The boundary conditions on the external boundary $\partial\Omega$, i.e. the second equation in (\ref{eq:interface_elasticity2}), will be imposed strongly. The remaining equations in (\ref{eq:interface_elasticity2}), i.e. the interface conditions on $\Gamma$, will be imposed by introduction of auxiliary variables on $\Omega_h^\Gamma$: the vector-valued $\mbf{p}$ (similar to the dual version of $\phi$-FEM for the Dirichlet boundary conditions above) and matrix-valued $\mbf{y}_1$,$\mbf{y}_2$ (similar to  $\phi$-FEM for the Neumann boundary conditions). This gives, cf. the last two equations in \eqref{eq:interface_elasticity2}: 
\begin{align}
\label{InterMix1}
\mbf{u}_1 - \mbf{u}_2 + \mbf{p}\phi = 0 \,, \quad &\text{ on } \ \Omega_h^{\Gamma}, \\
\label{InterMix2}
\mbf{y}_i + \mbf{\sigma}_i (\mbf{u}_i) = 0  \,, \quad &\text{ on } \ \Omega_h^{\Gamma},\ i=1,2, \\ 
 \label{InterMix3}
\mbf{y}_1 \nabla	 \phi-  \mbf{y}_2\nabla \phi = 0 \,, \quad &\text{ on } \ \Omega_h^{\Gamma} \,.
\end{align}
Equation (\ref{InterMix3}) above extends the last equation in (\ref{eq:interface_elasticity2}) from $\Gamma$ to $\Omega_h^\Gamma$ since the normal on $\Gamma$ is colinear with the vector $\nabla\phi$ there.\phantom{\ref{InterMix2}}


We are now going to discretize equations (\ref{weakExtomi})--(\ref{InterMix3}). We fix an integer $k\ge 1$ and introduce the FE spaces for the primary variables $\mbf{u}_i$: 
\begin{multline}\label{espaceVhi}
     V_{h,i} := \big\lbrace \mbf{v}_h:\Omega_{h,i}\to\mathbb{R}^d : \mbf{v}_{h |T} \in \mathbb{P}^k(T)^d \ \ \forall T \in \mathcal{T}_h, \ \mbf{v}_h\text{ continuous on }\Omega_{h,i}\,, \\ 
     \text{ and } \mbf{v}_h = I_h\mbf{u}^g\ \text{ on } \partial \Omega \big\rbrace \,
\end{multline}
with the standard FE interpolation $I_h$, and their homogeneous counterparts $V_{h,i}^0$ with the constraint $\mbf{v}_h = \mbf{0}\ \text{ on } \partial \Omega$, to be used for the test functions. We recall moreover the spaces $Q_{h}(\Omega_h^{\Gamma})$ and $Z_h(\Omega_h^{\Gamma})$ defined respectively by  (\ref{spaceQh}) and  (\ref{espaceZh}). Combining (\ref{weakExtomi}) with (\ref{InterMix1})--(\ref{InterMix3}) taken in the least square sense, gives the following scheme:\\
find $\mbf{u}_{h,1}\in V_{h,1}$, $\mbf{u}_{h,2}\in V_{h,2}$, $\mbf{p}_h \in Q_{h}^k(\Omega_h^{\Gamma})$, $\mbf{y}_{h,1}, \mbf{y}_{h,2} \in Z_h(\Omega_h^{\Gamma})$  such that, 
\begin{multline}\label{discreteInterf}
\sum_{i=1}^2 \int_{\Omega_{h,i}} \mbf{\sigma}_i(\mbf{u}_{h,i}):\nabla\mbf{v}_{h,i} 
+ \sum_{i=1}^2\int_{\partial\Omega_{h,i}}\mbf{y}_{h,i}\mbf{n} \cdot \mbf{v}_h \\
 + \frac{\gamma_p}{h^2} \int_{\Omega_h^{\Gamma}} (\mbf{u}_{h,1} - \mbf{u}_{h,2} + \frac{1}{h} \mbf{p}_h\phi_h)\cdot( \mbf{v}_{h,1} - \mbf{v}_{h,2} + \frac{1}{h} \mbf{q}_h\phi_h) \\
+ \gamma_u \sum_{i=1}^2\int_{\Omega_h^{\Gamma}}(\mbf{y}_{h,i}+ \mbf{\sigma}_i(\mbf{u}_{h,i})) : (\mbf{z}_{h,i}+ \mbf{\sigma}_i(\mbf{v}_{h,i})) \\ 
+ \frac{\gamma_y}{h^2}  \int_{\Omega_h^{\Gamma}} ( \mbf{y}_{h,1}\nabla \phi_h -  \mbf{y}_{h,2}\nabla \phi_h )\cdot( \mbf{z}_{h,1}\nabla \phi_h- \mbf{z}_{h,2} \nabla \phi_h) \\
+ \sum_{i=1}^2 \left(G_h(\mbf{u}_{h,i}, \mbf{v}_{h,i}) + J_h^{lhs,N} (\mbf{y}_{h,i}, \mbf{z}_{h,i}) \right)
 = \sum_{i=1}^2\int_{\Omega_{h,i}} \mbf{f} \cdot \mbf{v}_{h,i} 
 + \sum_{i=1}^2 J_h^{rhs,N} (\mbf{z}_{h,i}) \,, \\\quad \forall \, 
 \mbf{v}_{h,1}\in V_{h,1}^0, \mbf{v}_{h,2}\in V_{h,2}^0, \mbf{q}_h \in Q_{h}^k(\Omega_h^{\Gamma}), \mbf{z}_{h,1}, \mbf{z}_{h,2} \in Z_h(\Omega_h^{\Gamma}) \,.
\end{multline}
Similarly to the previous settings, we have added here the ghost stabilization $G_h$ defined by (\ref{Gh}) and the additional stabilization $J_h^{rhs,N}$ defined by (\ref{JhN}) with $\Omega_h^{\Gamma_N}$ replaced by $\Omega_h^{\Gamma}$ and imposing  $\Div\mbf{y}_i=\mbf{f}$ on  $\Omega_h^{\Gamma}$ in the least squares sense. 

\begin{figure}[btp]
\centering
\includegraphics[width=0.45\textwidth]{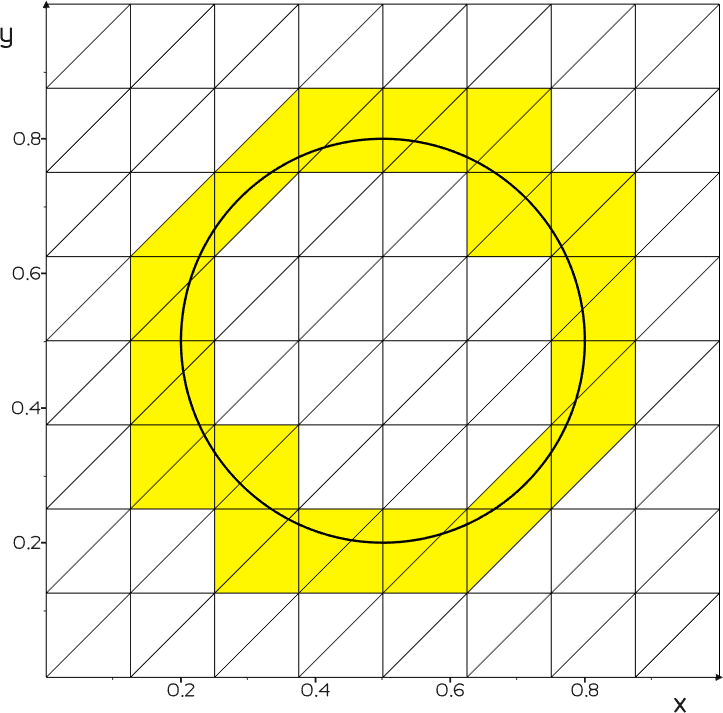}
\includegraphics[width=0.45\textwidth]{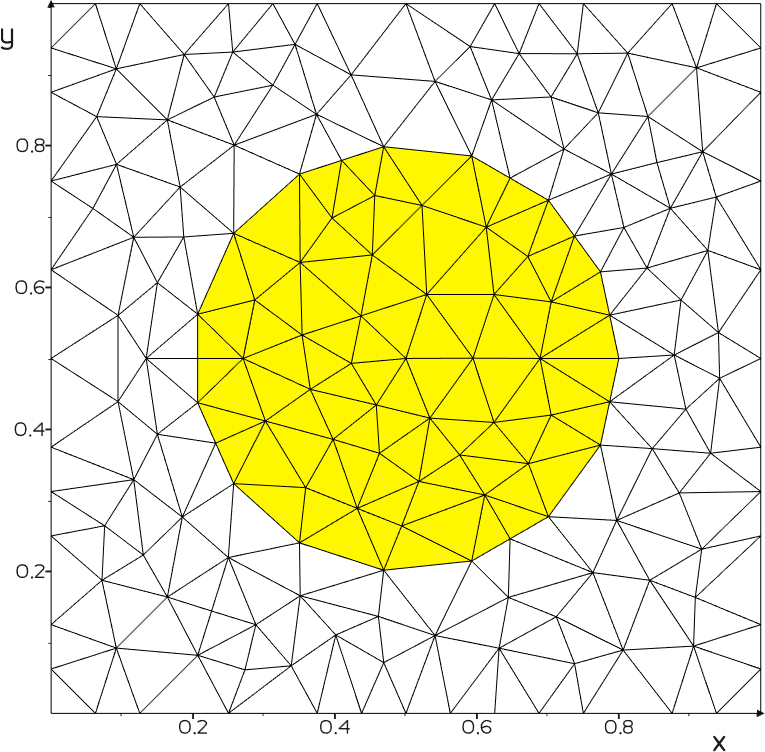}
\caption{Linear elasticity with multiple materials. Left: a mesh used for $\phi$-FEM ($\Omega_h^{\Gamma}$ painted in yellow); Right: a mesh matching the interface for standard FEM (yellow and white  represent the two materials).  }\label{fig:interface mesh}
\end{figure}
\begin{figure}[btp]
\centering
\begin{tikzpicture}
\begin{loglogaxis}[name = ax1, width = .45\textwidth, xlabel = $h$, ylabel = relative error, 
			legend style = { at={(1,1.4)}}]
\addplot coordinates{
(0.14142135623730964,5.675910662284101e-05)
(0.07071067811865482,6.2668613054719945e-06)
(0.03535533905932741,7.89051810290236e-07)
(0.01767766952966378,9.732876356767048e-08)
(0.00883883476483197,1.2751093753131592e-08)
};
\addplot coordinates{
(0.14142135623730964,0.0005631985529905538)
(0.07071067811865482,0.00013288249461713975)
(0.03535533905932741,3.2590636010544636e-05)
(0.01767766952966378,8.08676312960304e-06)
(0.00883883476483197,2.014362599383071e-06)};
\addplot coordinates{
(0.1552914270615126,0.011092067908400144)
(0.07940855385524764,0.002876516858202795)
(0.03977439957208735,0.0007593967526539493)
(0.019883151483778883,0.00018689260740527478)
(0.009943258702793148,4.682948959033794e-05)};
\addplot coordinates{
(0.1552914270615126,0.040116831985460386)
(0.07940855385524764,0.014226141479733757)
(0.03977439957208735,0.005225673219304294)
(0.019883151483778883,0.0018109814763434579)
(0.009943258702793148,0.0006366421849599564)};
\logLogSlopeTriangle{0.3}{0.2}{0.3}{2}{red};
\logLogSlopeTriangle{0.53}{0.2}{0.12}{3}{blue};
\legend{$L^2$ error $\phi$-FEM,$H^1$ error $\phi$-FEM,$L^2$ error standard FEM,$H^1$ error standard FEM}
\end{loglogaxis}
\begin{loglogaxis}[name = ax2, at = {(ax1.south east)}, xshift = 2cm, width = .45\textwidth, 
			ylabel = Computing time (s), xlabel = $L^2$ relative error,
			legend style = { at={(1,1.1)}, legend columns =2,
			/tikz/column 2/.style={column sep = 10pt}}]
\addlegendentry{Standard FEM}	
\addplot coordinates{
(0.03945192586664578,0.15355443954467773)
(0.011092067908400144,0.24408245086669922)
(0.002876516858202795,0.5149095058441162)
(0.0007594702182368449,1.0383923053741455)
(0.00018690079814396443,3.6394455432891846)
};
\addlegendentry{$\phi$-FEM}
\addplot coordinates{
(5.675910662284101e-05,0.26125335693359375)
(6.2668613054719945e-06,0.37996554374694824)
(7.89051810290236e-07,0.7896356582641602)
(9.732876356767048e-08,3.172820568084717)
(1.2751093753131592e-08,9.983308553695679)
};
\end{loglogaxis}
\end{tikzpicture}
\caption{Test case with multiple materials. Left: $H^1$ and $L^2$ relative error obtained with $\phi$-FEM and the standard FEM. Right: computing times for $\phi$-FEM and the standard FEM.}\label{fig:error_interface} 
\end{figure}
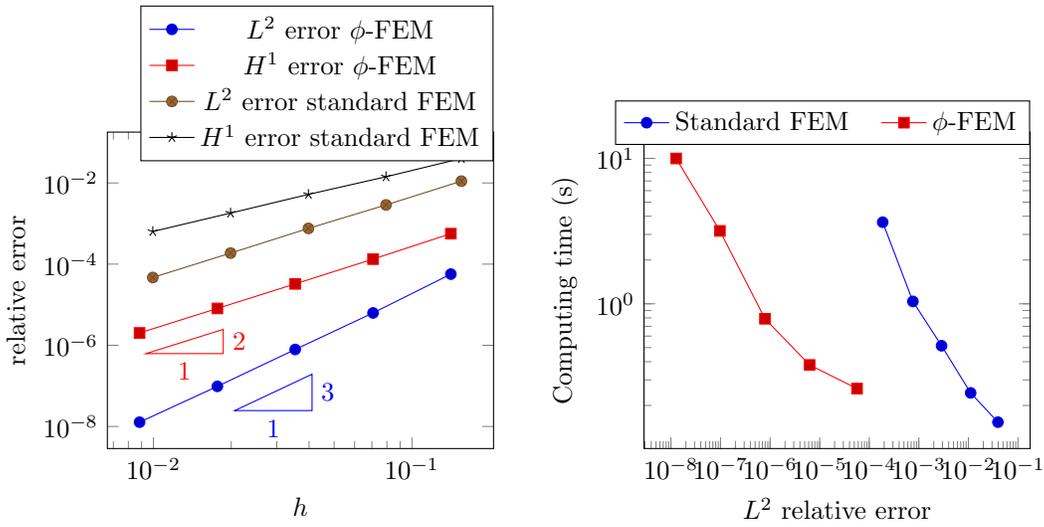

\paragraph*{Test case: }
Consider $\Omega=(0,1)^2$ and $\Omega_1$, $\Omega_2$ defined by the the level-set $\phi$  
\begin{equation*}
\phi(x,y) = - R^2 + (x-0.5)^2 + (y-0.5)^2\,,
\end{equation*}
with $R=0.3$ as illustrated on Fig.~\ref{fig:interface}. We want to solve \eqref{eq:interface_elasticity} with  the manufactured  radial solution
\[ 
\mbf{u}=\mbf{u}_{ex} = \left\{\begin{array}{ll}
\frac{1}{E_1}(\cos(r)-\cos(R))(1,1)^T&\mbox{ if }r<R,\\
\frac{1}{E_2}(\cos(r)-\cos(R))(1,1)^T&\mbox{ else,}
\end{array}\right.\,,
\] 
where $r=\sqrt{(x-0.5)^2+(y-0.5)^2}$.
Thus
\[\mbf{f}=-\Div(\sigma_1((\cos(r)-\cos(R))(1,1)^T)/E_1\]
and $\mbf{u}_g=\mbf{u}_{ex}$.

The material parameters are given by $E_1 = 7$, $E_2=2.28$ and $\nu_1 = \nu_2 = 0.3$. The meshes used for $\phi$-FEM and for the standard FEM are illustrated in Fig.~\ref{fig:interface mesh}. In the latter case, the mesh should resolve the interface $r=R$ so that the solution $\mbf{u}_h\in V_h$ is obtained by the straight-forward scheme
\begin{equation}\label{FEMinterf}
\sum_{i=1}^2 \int_{\Omega_{h,i}} \mbf{\sigma}_i(\mbf{u}_h):\nabla\mbf{v}_h
 = \int_{\Omega} \mbf{f} \cdot \mbf{v}_h,~\forall~ \mbf{v}_h\in V_h^0\,,
\end{equation}
where $V_h$ is the conforming $\mathbb{P}^k$ FE space approximating $\mbf{u}_g$ on $\partial\Omega$ and  $V_h^0$ is its homogeneous analogue. The results obtained with $\phi$-FEM (\ref{discreteInterf}) and FEM (\ref{FEMinterf}) using $\mathbb{P}^2$ piecewise polynomials ($k=2$) are reported in Fig.~\ref{fig:error_interface}. The conclusions remain the same as in the previous setting: $\phi$-FEM is  more precise on comparable meshes and less expensive in terms of the computing times for a given error tolerance. 

\section{Linear elasticity with cracks }\label{sectFracture}
We now want to consider the linear elasticity problem posed on a cracked domain $\Omega\setminus\Gamma_f$ with $\Gamma_f$ being a line (a surface) inside $\Omega$:
\begin{equation}\label{eq:fracture_elasticity0}
\begin{cases}
- \Div \mbf{\sigma} (\mbf{u}) &= \mbf{f} \,, \text{ on } \ \Omega\setminus\Gamma_f \,, \\
\mbf{u} &= \mbf{u}^g \,, \text{ on } \ \partial \Omega \,, \\
\mbf{\sigma}(\mbf{u})\mbf{n} &= \mbf{g} \,, \text{ on } \ \Gamma_{f} \,.
\end{cases}
\end{equation}
This problem is actually what XFEM was originally designed for, cf. \cite{moes}. We are now going to adapt $\phi$-FEM to it.

\begin{figure}[b]
\centering
\begin{tikzpicture}
\begin{axis}[
	xmin = 0, xmax = 1,
	ymin = 0, ymax = 1.0,  xtick={0,1},ytick={1}, width = 0.45\textwidth]
	\addplot[ domain = 0:0.5,samples = 1000, line width=1pt, cyan] { 0.25*( sin(deg(2*pi*x))) + 0.5}
	node [pos=0.5, above left]{\Large$\Gamma_{int}$};
	\addplot[ domain = 0.5:1,samples = 1000, line width=1pt, red] { 0.25*( sin(deg(2*pi*x))) + 0.5}
	node[pos = 0.7, above left]{\Large$\Gamma_f$};

	\draw(0.7,0.7)node[black]{\Large$\Omega_1$};
	\draw(0.3,0.3)node[black]{\Large$\Omega_2$};
		\end{axis}

	\draw(2.5,-0.1)node[below]{$\mbf{u} = \mbf{u}^g$};

\end{tikzpicture}
\caption{Geometry notations to represent the crack. $\Gamma_{int}$  and $\Gamma_f$ represents the fictitious interface and the actual crack respectively.}\label{fig:figurefracture_phifem}
\end{figure}
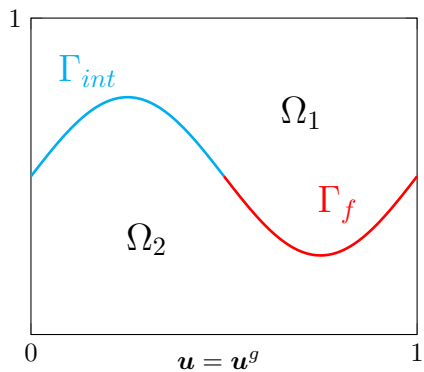

In practice, the crack geometry is given by the primary level set $\phi$ (to locate the line or surface of the crack) and the secondary level set $\psi$ (to locate the tip or the front of the crack):
$$
\Gamma_{f} := \Omega \cap  \lbrace \phi=0 \rbrace \cap \lbrace \psi < 0 \rbrace \,.
$$
To fix the ideas, let us suppose that the line (surface) $\Gamma:=\{\phi=0\}$ splits $\Omega$ into two sub-domains  $\Omega_1$ and $\Omega_2$, characterized by $\{\phi<0\}$ and $\{\phi>0\}$ respectively, as illustrated at Fig.~\ref{fig:figurefracture_phifem}. The interface $\Gamma$ thus consists of the fracture location $\Gamma_f$ and the remaining  (fictitious) part $ \Gamma_{int}$: 
\[ \Gamma_{int} := \Omega \cap  \lbrace \phi=0 \rbrace \cap \lbrace \psi > 0 \rbrace \,.
\]
In order to reuse the $\phi$-FEM scheme (\ref{discreteInterf}) introduced for the interface problem above, we reformulate problem (\ref{eq:fracture_elasticity0}) in terms of two separate unknowns $\mbf{u}_i=\mbf{u} |_{\Omega_i}$, $i=1,2$:
\begin{equation}\label{eq:fracture_elasticity}
\begin{cases}
- \Div \mbf{\sigma} (\mbf{u}_i) &= \mbf{f} \,, \text{ on } \ \Omega_i \,, \\
\mbf{u}_i &= \mbf{u}^g \,, \text{ on } \ \partial \Omega \,, \\
[\mbf{u} ] &= 0\,, \text{ on } \ \Gamma_{int} \,, \\
[\mbf{\sigma}(\mbf{u})\mbf{n}] &= 0 \,, \text{ on } \ \Gamma_{int} \,, \\
\mbf{\sigma}(\mbf{u})\mbf{n} &= \mbf{g} \,, \text{ on } \ \Gamma_{f} \,.
\end{cases}
\end{equation}

\begin{figure}[btp]
\centering
\includegraphics[width=0.45\textwidth]{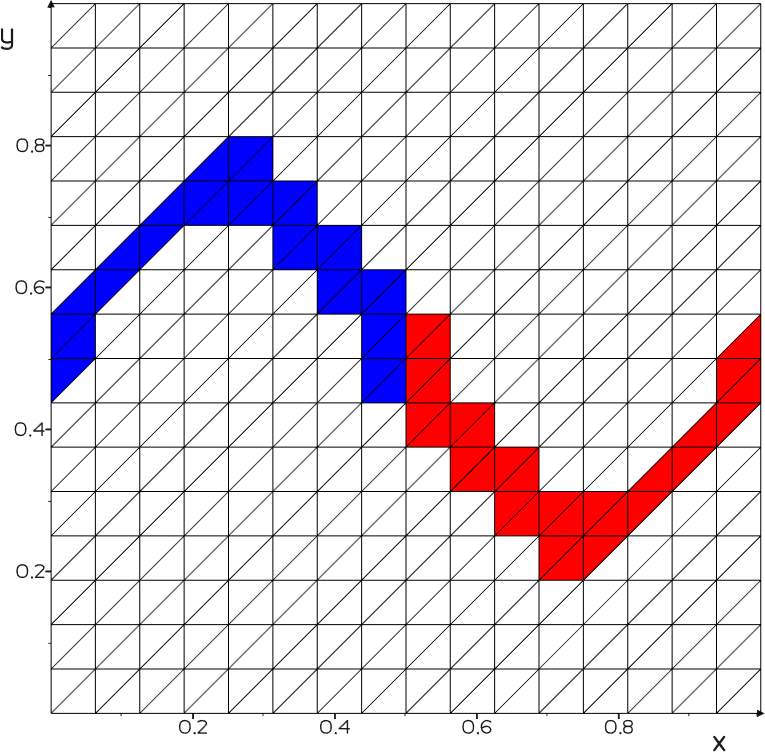}
\includegraphics[width=0.45\textwidth]{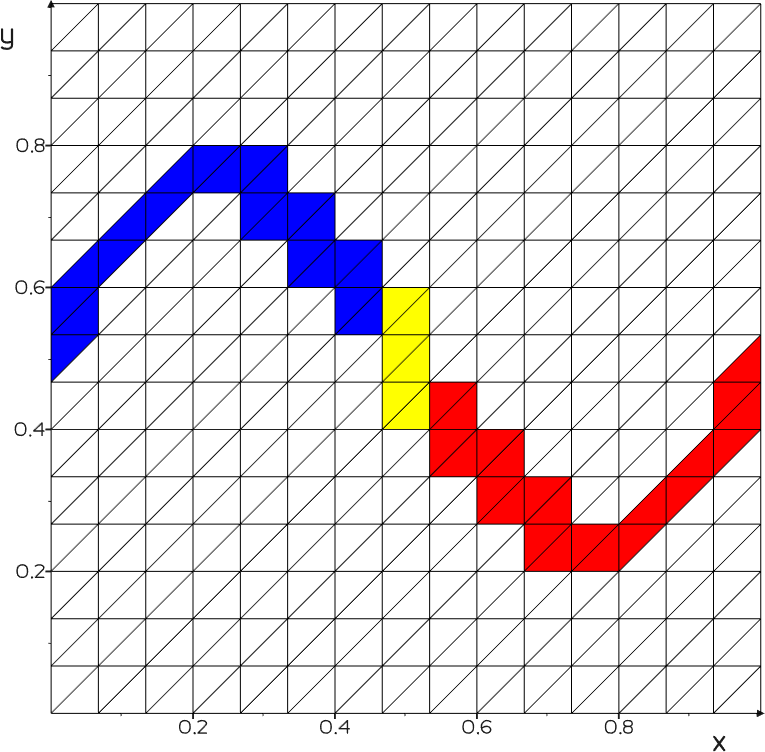}
\caption{Test case with a crack, meshes used for $\phi$-FEM.  Left: a mesh resolving the crack tip; the cells in $\mathcal{T}_{h}^{\Gamma_{int}}$ in blue; the cells in $\mathcal{T}_{h}^{\Gamma_{f}}$ in red. Right: a mesh not  resolving the crack tip; in addition to blue and red cells, there are yellow cells not belonging to $\mathcal{T}_{h}^{\Gamma_{int}}$ or $\mathcal{T}_{h}^{\Gamma_{f}}$.  }\label{fig:meshes fracture}
\end{figure}

We are interested again in a situation where $\Omega$ is sufficiently simple-shaped so that a matching mesh $T_h$ on $\Omega$ is easily available, but this mesh does not  match the internal interface $\Gamma$. As in the preceding section, we are thus going to discretize separately $\mbf{u}_1$  on $\Omega_1$ and $\mbf{u}_2$ on $\Omega_2$ starting from the reformulation (\ref{eq:fracture_elasticity}).  To this end, we introduce two active sub-meshes $\mathcal{T}_{h,1}$, $\mathcal{T}_{h,2}$ as in (\ref{Th12}), based on the piecewise polynomial approximation $\phi_h$ of $\phi$. 
 We also introduce the interface mesh $\mathcal{T}_{h}^\Gamma=\mathcal{T}_{h,1}\cap\mathcal{T}_{h,2}$\,, which we further split into two sub-meshes with respect to the secondary level set $\psi$, similarly to our treatment of the mixed boundary conditions, cf. (\ref{ThGamDN}):
\begin{equation*}
    \mathcal{T}_h^{\Gamma_f} := \{ T \in \mathcal{T}_h^{\Gamma} : \psi\leqslant 0\text{ on }T \} \qquad \text{ and } \qquad 
\mathcal{T}_h^{\Gamma_{int}} := \{ T \in \mathcal{T}_h^{\Gamma} : \psi\geqslant 0\text{ on }T \}\,.
\end{equation*} 
Note that there may be some cells in $\mathcal{T}_h^{\Gamma}$ that are not marked as either $\mathcal{T}_h^{\Gamma_f}$ or $\mathcal{T}_h^{\Gamma_{int}}$. This is illustrated by the mesh example on the right of Fig.~\ref{fig:meshes fracture}, where the cells in $\mathcal{T}_h^{\Gamma_f}$ and $\mathcal{T}_h^{\Gamma_{int}}$ are painted in red and blue respectively, but there remain some cells  $\mathcal{T}_h^{\Gamma}$ that are in neither of these categories. The are painted in yellow on the picture. These are the cells intersected by the line $\{\psi=0\}$. The crack tip happens to be thus inside one of the yellow cells. 

Everything is now set up to adapt the $\phi$-FEM approaches of the two preceding sections to the equations (\ref{eq:fracture_elasticity}). We choose an integer $k\ge 1$ and introduce first the FE spaces $V_{h,1}$, $V_{h,2}$ together with their homogeneous counterparts $V_{h,1}^0$, $V_{h,2}^0$ as in (\ref{espaceVhi}) to approximate $\mbf{u}_1$ and $\mbf{u}_2$. These will be used in the discretization of the variational formulation of the first equation in (\ref{eq:fracture_elasticity}) together with the boundary conditions on $\partial\Omega$. The remaining equations in (\ref{eq:fracture_elasticity}), i.e. the relations on $\Gamma_{int}$ and $\Gamma_f$ will be treated by the introduction of auxiliary variables on the appropriate parts of $\Omega_h^\Gamma$ (the domain of the mesh $\Th^\Gamma$):
\begin{itemize}
    \item the vector-valued unknown $\mbf{p}$ and the matrix-valued unknowns $\mbf{y}_1$,$\mbf{y}_2$ on $\Omega_h^{\Gamma_{int}}$ (the domain of the mesh $\Th^{\Gamma_{int}}$). These will serve to impose the continuity of both the displacement and the normal force on ${\Gamma_{int}}$ thorough the equations
    \begin{align*}
  \mbf{u}_1 - \mbf{u}_2 + \mbf{p}\phi = 0 \,, \quad &\text{ on } \ \Omega_h^{\Gamma_{int}}\,, \\
  \mbf{y}_i = - \mbf{\sigma} (\mbf{u}_i)  \,, \quad &\text{ on } \ \Omega_h^{\Gamma_{int}}\,, \\ 
  \mbf{y}_1 \cdot \nabla	 \phi-  \mbf{y}_2 \cdot \nabla \phi = 0 \,, \quad &\text{ on } \ \Omega_h^{\Gamma_{int}} \,,
     \end{align*}
    which are exactly the same as (\ref{InterMix1})--(\ref{InterMix3}) with the only exception that they are posed on the appropriate portion of $\Omega_h^\Gamma$ rather than on entire $\Omega_h^\Gamma$. These variables will be discretized in FE spaces $Q_{h}^k(\Omega_h^{\Gamma_{int}})$ for $\mbf{p}$ and $Z_h(\Omega_h^{\Gamma_{int}})$ for $\mbf{y}_1$,$\mbf{y}_2$, defined by (\ref{spaceQh}) and (\ref{espaceZh}) respectively. 
    
    \item the vector-valued unknowns $\mbf{p}_i^N$ and the matrix-valued unknown $\mbf{y}_i^N$, $i=1,2$ on  $\Omega_h^{\Gamma_f}$ (the domain of the mesh $\Th^{\Gamma_f}$). These will serve to impose the Neumann boundary conditions on both sides of ${\Gamma_f}$ thorough the equations
    \begin{align*}
  \mbf{y}_i^N = - \mbf{\sigma} (\mbf{u}_i)  \,, \quad &\text{ on } \ \Omega_h^{\Gamma_{f}}\,, \\ 
  \mbf{y}_i^N \nabla \phi + \mbf{p}_i^N \phi + \mbf{g} |\nabla \phi | = 0 \,, \quad &\text{ on } \ \Omega_h^{\Gamma_{f}} \,,
    \end{align*}  
    which are exactly the same as (\ref{phiN}a-b) with the only exception that the domain is renamed to  $\Omega_h^{\Gamma_f}$ from $\Omega_h^{\Gamma_N}$. These variables will be discretized in FE spaces $Q_{h}^{k-1}(\Omega_h^{\Gamma_f})$ for $\mbf{p}_i^N$ and $Z_{h}(\Omega_h^{\Gamma_f})$ for $\mbf{y}_i^N$, defined again by (\ref{spaceQh}) and (\ref{espaceZh}) respectively.   
\end{itemize}
Note that the combination of equations above does not impose the appropriate interface conditions on the whole of $\Gamma$ since the latter may be not completely covered by  $\Omega_h^{\Gamma_f}\cup\Omega_h^{\Gamma_{int}}$. Fortunately, this defect of the formulation on the continuous level can be repaired on the discrete level by adding the appropriate stabilization to the FE discretization, similarly to what we have already seen in the setting with mixed boundary conditions.

All this results in the following FE scheme:
find $\mbf{u}_{h,1}\in V_{h,1}$, $\mbf{u}_{h,2}\in V_{h,2}$, $\mbf{p}_h\in Q_{h}^{k}(\Omega_h^{\Gamma_D})$, $\mbf{y}_{h,1}, \mbf{y}_{h,2}\in Z_{h}(\Omega_h^{\Gamma_N})$, $\mbf{p}_{h,1}^N, \mbf{p}_{h,2}^N \in Q_{h}^{k-1}(\Omega_h^{\Gamma_N})$, $\mbf{y}_{h,1}^N,\mbf{y}_{h,2}^N\in Z_{h}(\Omega_h^{\Gamma_N})$ such that
\begin{multline} 
\sum_{i = 1}^2 \big( \int_{\Omega_{h, i}} \mbf{\sigma} (\mbf{u}_{h, i}) : \mbf{\nabla} \mbf{v}_{h, i} + \int_{\partial \Omega_{h, i, int}} \mbf{y}_{h, i} \mbf{n} \cdot \mbf{v}_{h, i} + \int_{\partial \Omega_{h, i, f}} \mbf{y}_{h, i}^N \mbf{n} \cdot \mbf{v}_{h, i} \\
- \int_{\partial\Omega_{h, i}\setminus(\partial\Omega_{h, i, int}\cup\partial\Omega_{h, i, f})} \mbf{\sigma} (\mbf{u}_{h, i})\mbf{n} \cdot \mbf{v}_{h, i} \big) \\ 
+ \frac{\gamma_p}{h^2}  \int_{\Omega_h^{\Gamma_{int}}} (\mbf{u}_{h, 1} - \mbf{u}_{h, 2} + \frac{1}{h} \mbf{p}_h \phi_h) \cdot (\mbf{v}_{h, 1} - \mbf{v}_{h, 2} + \frac{1}{h} \mbf{q}_h \phi_h) \\ 
+ \gamma_u  \sum_{i = 1}^2 \int_{\Omega_h^{\Gamma_{int}}} (\mbf{y}_{h, i} + \mbf{\sigma} (\mbf{u}_{h, i})) : (\mbf{z}_{h, i} + \mbf{\sigma} (\mbf{v}_{h, i})) \\
+ \frac{\gamma_y}{h^2}  \int_{\Omega_h^{\Gamma_{int}}} ( \mbf{y}_{h,1}\nabla \phi_h -  \mbf{y}_{h,2}\nabla \phi_h )\cdot( \mbf{z}_{h,1}\nabla \phi_h- \mbf{z}_{h,2} \nabla \phi_h) 
\\
+ \gamma_{u, N}  \sum_{i = 1}^2 \int_{\Omega_h^{\Gamma_f}} (\mbf{y}_{h, i}^N + \mbf{\sigma} (\mbf{u}_{h, i})) : (\mbf{z}_{h, i}^N + \mbf{\sigma} (\mbf{v}_{h, i})) \\ + \frac{\gamma_{p, N}}{h^2}  \sum_{i = 1}^2 \int_{\Omega_h^{\Gamma_f}} (\mbf{y}_{h, i}^N \nabla \phi_h + \frac{1}{h} \mbf{p}_{h, i}^N \phi_h) \cdot (\mbf{z}_{h, i}^N \nabla \phi_h + \frac{1}{h} \mbf{q}_{h, i}^N \phi_h) \\ 
+ \sum_{i = 1}^2 \left( G_h  \left( \mbf{u}_{h, i}, \mbf{v}_{h, i} \right) 
 + J_h^{{lhs}, {int}} \left( \mbf{y}_{h, i}, \mbf{z}_{h, i} \right) 
 + J_h^{{lhs}, f} \left( \mbf{y}_{h, i}^N, {\mbf{z }_{h, i}^N}  \right) \right) \\ 
 = \sum_{i = 1}^2 \int_{\Omega_{h, i}} \mbf{f} \cdot \mbf{v}_{h, i} - \frac{\gamma_{p, N}}{h^2}  \sum_{i = 1}^2 \int_{\Omega_h^{\Gamma_f}} \mbf{g} | \nabla \phi_h |  (\mbf{z}_{h, i}^N \nabla \phi_h + \frac{1}{h} \mbf{q}_{h, i}^N \phi_h) \\ 
+ \sum_{i = 1}^2 \left( 
J_h^{{rhs}, {int}} \left( \mbf{z}_{h, i} \right) 
 + J_h^{{rhs}, f} \left( {\mbf{z }_{h, i}^N}  \right) \right), \\
\forall \mbf{v}_{h, 1} \in V_{h, 1}^0, \mbf{v}_{h, 2} \in V_{h, 2}^0, \mbf{q}_h \in Q_{h}^{k}(\Omega_h^{\Gamma_D}), \mbf{z}_{h, 1}, \mbf{z}_{h, 2} \in Z_h, \mbf{q}_{h, 1}^N, \mbf{q}_{h, 2}^N \in Q_{h}^{k-1}(\Omega_h^{\Gamma_N}),\\ \mbf{\mbf{z}_{h, 1}^N, z}_{h, 2}^N \in Z_{h}(\Omega_h^{\Gamma_N}) \,. 
\label{discreteCrack}
\end{multline}
As usual,  we have added here the ghost stabilization $G_h$ (\ref{Gh}) and the additional stabilizations  $J_h^{lhs,int}$, $J_h^{lhs,f}$ 
 (accompanied by the their counterparts on the right-hand side for the consistency) that are copied from the $J_h^{lhs,N}$ in (\ref{JhN}) but adjusted to the corresponding sub-meshes:
$$    J_h^{lhs,int}(\mbf{y},\mbf{z}) =\gamma_{div} \int_{\Omega_h^{\Gamma_{int}}}\Div \mbf{y} \cdot \Div \mbf{z}
\,,
\quad J_h^{lhs,f}(\mbf{y},\mbf{z}) =\gamma_{div} \int_{\Omega_h^{\Gamma_f}}\Div \mbf{y} \cdot \Div \mbf{z}\,.
$$
The boundary integrals are rewritten  in terms of $\mbf{y}_i$,$\mbf{y}_i^N$ wherever possible. We have here denoted by $\partial\Omega_{h,i}$ the part of the boundary of $\Omega_{h,i}$ other than $\partial\Omega$ and introduced $\partial\Omega_{h,i,int}$ as the part of $\partial\Omega_{h,i}$ formed by the boundary facets of $\mathcal{T}_{h,i}$  belonging to the cells in $\mathcal{T}_{h}^{\Gamma_{int}}$. The same for $\partial\Omega_{h,i,f}$.

\begin{figure}[btp]
\centering
\begin{tikzpicture}
\begin{loglogaxis}[name = plot, width = .45\textwidth, ymin=3e-8, ymax=5e-3, xlabel = $h$, ylabel = Relative error,legend pos=north west]
\addplot coordinates{
(0.07071067811865482,9.777452793013566e-05)
(0.03535533905932741,1.0419211256548806e-05)
(0.01767766952966378,2.967883868554893e-06)
(0.00883883476483197,4.3827759106063176e-07)
(0.004419417382415985,5.625169839168734e-08)};
\addplot coordinates{
(0.07071067811865482,0.0006052966265906567)
(0.03535533905932741,0.00011243253507143713)
(0.01767766952966378,3.6174590887733385e-05)
(0.00883883476483197,5.723173873151866e-06)
(0.004419417382415985,8.466975252389401e-07)};
\logLogSlopeTriangle{0.53}{0.2}{0.46}{2}{red};
\logLogSlopeTriangle{0.53}{0.2}{0.23}{3}{blue};
\legend{$L^2$ error, $H^1$ error}
\end{loglogaxis}
\end{tikzpicture}
\begin{tikzpicture}
\begin{loglogaxis}[name = plot, width = .45\textwidth, ymin=3e-8, ymax=5e-3, xlabel = $h$, ylabel = Relative error,legend pos=north west]
\addplot coordinates{
(0.06734350297014746,0.00011415468313107934)
(0.034493013716416984,1.7495400368412347e-05)
(0.017459426695964213,2.954050638336371e-06)
(0.008783935170019241,4.182817999120714e-07)
(0.004405649727018988,5.5993729055912355e-08)
};
\addplot coordinates{
(0.06734350297014746,0.0006725963282258503)
(0.034493013716416984,0.00013063709444329212)
(0.017459426695964213,3.121288024010827e-05)
(0.008783935170019241,5.393354911850969e-06)
(0.004405649727018988,8.035075055809038e-07)
};
\logLogSlopeTriangle{0.53}{0.2}{0.46}{2}{red};
\logLogSlopeTriangle{0.53}{0.2}{0.23}{3}{blue};
\legend{$L^2$ error, $H^1$ error}
\end{loglogaxis}
\end{tikzpicture}
\caption{Test case with a crack, $H^1$ and $L^2$ relative errors. Left: on meshes resolving the crack tip.  Right: on meshes not resolving the crack tip. }\label{fig:convergence_fracture}
\end{figure}
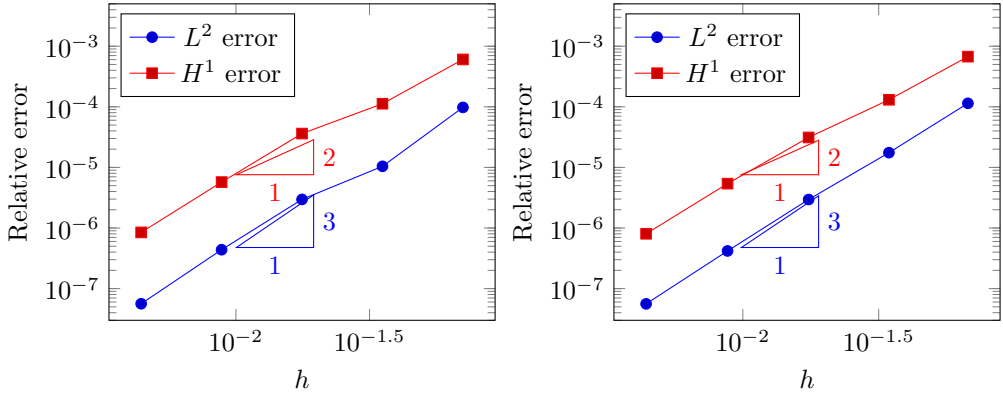

\paragraph*{Test case: }
Let $\Omega = (0,1)^2$ and the interface $\Gamma$ be given by the level set 
\begin{equation*}
    \phi (x,y) = y - \frac{1}{4} \sin(2 \pi x) - \frac{1}{2}\,.
\end{equation*}
We choose the crack tip to be at $x=0.5$ so that
\[ \Gamma_{int} := \{\phi=0\}\cap \lbrace x < 0.5 \rbrace \  \quad \text{ and } \quad \Gamma_{f} := \{\phi=0\}\ \cap \lbrace x > 0.5 \rbrace \,. \]
This is the setting represented at Fig.~\ref{fig:figurefracture_phifem}.

We use the $\phi$-FEM (\ref{discreteCrack}) to solve \eqref{eq:fracture_elasticity0} with the manufactured solution
\[ \mbf{u} = \mbf{u}_{ex} = (\sin(x) \times \exp(y), \sin(y) \times \exp(x))^T \] 
which gives $\mbf{f}$, $\mbf{g}$, and $\mbf{u}^g$ by substitution. The force on the crack $\mbf{g}$ should be extended to a vicinity of $\Gamma_f$ and we implement it by
 \[
\mbf{g} = \mbf{\sigma}(\mbf{u}_{ex}) \frac{\nabla \phi}{\| \nabla \phi\|} + \phi\mbf{u}_{ex}  \,. \]
We choose $\gamma_u = \gamma_p = \gamma_{div} = \gamma_{u,N} = \gamma_{p,N} = \gamma_{div,N} = 1.0$, $\sigma_p = 1.0$ and $\sigma_D = 20.0$.

We have conducted two series of numerical experiments using $\phi$-FEM (\ref{discreteCrack}) with $\mathbb{P}^2$ Lagrange polynomials ($k=2$) on families of meshes presented at Fig.~\ref{fig:meshes fracture}, either resolving the crack tip (the mesh on the left) or not (the mesh on the right). The results are reported on Fig.~\ref{fig:convergence_fracture}. We see that $\phi$-FEM converges optimally, giving very similar results on both types of meshes.

\section{Heat equation}\label{sectHeat}

We finally demonstrate the applicability of the $\phi$-FEM approach to time-dependent problems. We take the example of the heat equation with Dirichlet boundary conditions: given a bounded domain $\Omega\in\mathbb{R}^d$, the initial conditions $u^0$ on $\Omega$, and the final time $T>0$, find the scalar field $u=u(x,t)$ such that
\begin{equation}\label{eq:heat}
\left\{\begin{array}{ll}
u_t - \Delta u = f& \mbox{ in } \Omega\times(0,T), \\
u=0&\mbox{ on }\Gamma\times(0,T), \\
u(.,0)=u^0&\mbox{ in }\Omega.
\end{array}\right.
\end{equation}

We are  interested again in the situation where a fitting mesh of $\Omega$ is not available. We rather assume that $\Omega$ is inscribed in a box ${\mathcal{O}}$ which is covered by a simple background mesh $\Th^{\mathcal{O}}$, and introduce the active mesh $\Th$ as in (\ref{Th}). We then follow the Direct Dirichlet $\phi$-FEM approach (\ref{DirectPhi}), (\ref{DirectDiscrete})
with the following modifications: 
\begin{itemize}
    \item We introduce the uniform partition of the time interval $I=[0,T]$ into time steps of length $\Delta t$ by the nodes $t_i=i\Delta t$. We discretize then \eqref{eq:heat} in time using implicit Euler scheme. On the continuous level this is formally written as: find $u^n$ (the approximation to $u$ at time $t_n$) in the form $u^n=\phi w^n$ successively for $n=1,2,\ldots$ solving
    \begin{equation}\label{DiscTime}
        \frac{\phi w^n - \phi w^{n-1}}{\Delta t} -\Delta(\phi w^n)=f^n
    \end{equation}
    where $f^n(\cdot)=f(t_n,\cdot)$.
    \item We extend (\ref{DiscTime}) to $\Omega_h$, integrate by parts on $\Omega_h$, and discretize the resulting variational formulation using a FE space and adding appropriate stabilizations. 
\end{itemize}

The $\phi$-FEM for \eqref{eq:heat} reads thus as:  find $w_h^n \in V_h$ for $n=1,2,\ldots$ with $V_h$ defined by (\ref{spaceVh}) such that 
\begin{multline}\label{discreteHet}
  \int_{\Omega_h} \frac{\phi_h w_h^n }{\Delta t} \phi_h v_h +
   \int_{\Omega_h} \nabla (\phi_h w_h^n) \cdot \nabla (\phi_h v_h) 
   - \int_{\partial \Omega_h} \frac{\partial}{\partial n}  (\phi_h w_h^{n}) \phi_h v_h \\
   + \sigma_D h \sum_{E \in \mathcal{F}_h^{\Gamma}} \int_E \left[ \partial_n(\phi_h w_h^n)\right] \cdot \left[ \partial_n(\phi_h v_h)\right] 
 - \sigma h^2 \sum_{T\in\Th^{\Gamma}}\int_{T} \left( \frac{\phi_h w_h^n
   }{\Delta t} - \Delta (\phi_h w_h^n) \right) \Delta (\phi_h v_h ) \\
 = \int_{\Omega_h} \left( \frac{\phi_h w_h^{n-1} }{\Delta t} + f^n \right)
   \phi_h v_h - \sigma h^2 \sum_{T\in\Th^{\Gamma}}\int_{T} \left( \frac{\phi_h w_h^{n-1}
   }{\Delta t} + f^n \right) \Delta (\phi_h v_h ) .
\end{multline}
We have added here the ghost stabilisation, similar to  (\ref{Gh}) but in simpler scalar setting, and additional stabilization inspired by  (\ref{Jhlrhs}). The idea for the latter is to take the governing equation in the strong form, which is now (\ref{DiscTime}), and to impose it in a least squares manner cell by cell. 


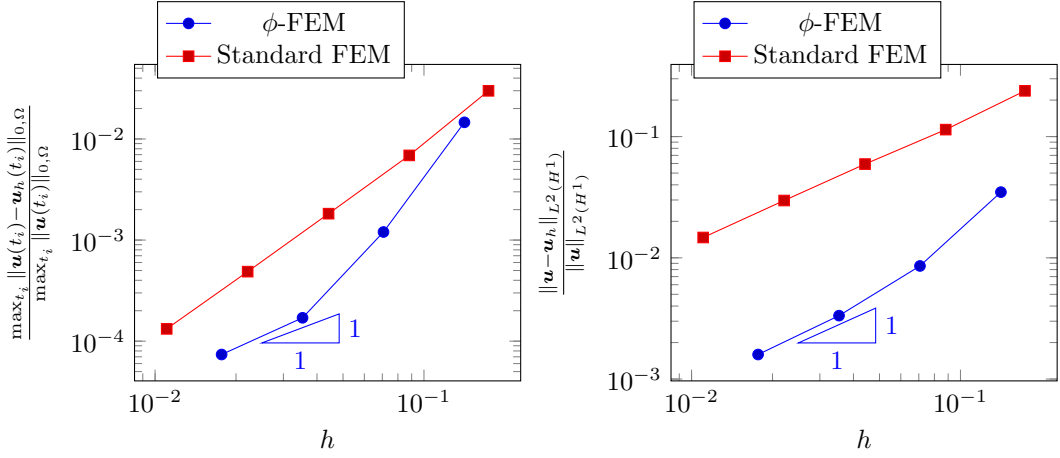
\begin{figure}[htbp]
\centering
\begin{tikzpicture}
\begin{loglogaxis}[name = ax1, width = .45\textwidth, xlabel = $h$, 
            ylabel = $\frac{\max_{t_i}\| \mbf{u}(t_i)-\mbf{u}_h(t_i) \|_{0,\Omega}}{\max_{t_i}\|\mbf{u}(t_i)\|_{0,\Omega}}$,
            legend style = { at={(0.7,1.2)}, legend columns =1,
			/tikz/column 2/.style={column sep = 10pt}}]
\addplot coordinates {
(0.14142135623730964,0.014613260310956685)
(0.07071067811865482,0.0012018592924686047)
(0.03535533905932741,0.00017023448816856176)
(0.01767766952966378,7.376483560822227e-05)
 };
\addplot coordinates { 
(0.1737929702933656,0.029956607570528575)
(0.08814237266741118,0.006878808462825571)
(0.044176082108213596,0.0018242512029954185)
(0.022095901417288732,0.0004871954091249405)
(0.01104624219043868,0.00013234126714823852)
};
\logLogSlopeTriangle{0.53}{0.2}{0.12}{1}{blue};
\legend{$\phi$-FEM, Standard FEM}
\end{loglogaxis}
\end{tikzpicture}
\begin{tikzpicture}
\begin{loglogaxis}[name = ax1, width = .45\textwidth, xlabel = $h$, 
            ylabel = $\frac{\| \mbf{u}-\mbf{u}_h \|_{L^2(H^1)}}{\|\mbf{u}\|_{L^2(H^1)}}$,
            legend style = { at={(0.7,1.2)}, legend columns =1,
			/tikz/column 2/.style={column sep = 10pt}}]

\addplot coordinates { 
(0.14142135623730964,0.03483439404608014)
(0.07071067811865482,0.008569471429621481)
(0.03535533905932741,0.0033349304158489697)
(0.01767766952966378,0.0015937182371751916)
};

\addplot coordinates{
(0.1737929702933656,0.23856788662022907)
(0.08814237266741118,0.11417978400648772)
(0.044176082108213596,0.05933004823534846)
(0.022095901417288732,0.02969997580122884)
(0.01104624219043868,0.0147143156195001)};

\logLogSlopeTriangle{0.53}{0.2}{0.12}{1}{blue};
\legend{ $\phi$-FEM, Standard FEM}
\end{loglogaxis}
\end{tikzpicture}
\caption{Test case for the heat equation; $\Delta t =h$. Left: $L^{\infty}(0,T;L^2(\Omega))$  relative errors.  Right: $L^2(0,T;H^1(\Omega))$ relative errors.}\label{fig:heat}
\end{figure}

\begin{figure}[htbp]
\centering
\begin{tikzpicture}
\begin{loglogaxis}[name = ax1, width = .45\textwidth, xlabel = $h$, 
            ylabel = $\frac{\max_{t_i}\| \mbf{u}(t_i)-\mbf{u}_h(t_i) \|_{0,\Omega}}{\max_{t_i}\|\mbf{u}(t_i)\|_{0,\Omega}}$,
            legend style = { at={(0.7,1.2)}, legend columns =1,
			/tikz/column 2/.style={column sep = 10pt}}]
\addplot coordinates {
(0.141421356237,0.0145429360786)
(0.0707106781187,0.00109257635152)
(0.0353553390593,0.000120961344644)
(0.0176776695297,3.0745030175e-05)
(0.00883883476483,7.81655564878e-06)
 };
\addplot coordinates { 
(0.141421356237,0.0347422759486)
(0.0707106781187,0.00840243388952)
(0.0353553390593,0.00333272225131)
(0.0176776695297,0.00159026149511)
(0.00883883476483,0.000788380527752)
};
\logLogSlopeTriangle{0.53}{0.2}{0.12}{2}{blue};
\legend{$\phi$-FEM, Standard FEM}
\end{loglogaxis}
\end{tikzpicture}
\begin{tikzpicture}
\begin{loglogaxis}[name = ax1, width = .45\textwidth, xlabel = $h$, 
            ylabel = $\frac{\| \mbf{u}-\mbf{u}_h \|_{L^2(H^1)}}{\|\mbf{u}\|_{L^2(H^1)}}$,
            legend style = { at={(0.7,1.2)}, legend columns =1,
			/tikz/column 2/.style={column sep = 10pt}}]

\addplot coordinates { 
(0.173792970293,0.0301514946318)
(0.0881423726674,0.00686051308391)
(0.0441760821082,0.00178094355391)
(0.0220959014173,0.000454820848822)
(0.0110462421904,0.000112234933215)
};

\addplot coordinates{
(0.173792970293,0.238572310108)
(0.0881423726674,0.114179466332)
(0.0441760821082,0.0593294708849)
(0.0220959014173,0.0296996298576)
(0.0110462421904,0.0147141332572)
};

\logLogSlopeTriangle{0.53}{0.2}{0.12}{2}{blue};
\legend{ $\phi$-FEM, Standard FEM}
\end{loglogaxis}
\end{tikzpicture}
\caption{Test case for the heat equation; $\Delta t =10 h^2$. Left: $L^{\infty}(0,T;L^2(\Omega))$  relative errors.  Right: $L^2(0,T;H^1(\Omega))$ relative errors.}\label{fig:heat2}
\end{figure}

\paragraph*{Test case: } We consider again the geometry of $\Omega$ and of the surrounding box $\mathcal{O}$ as in our first test case on page~\pageref{testcaseDir}. In particular, the level set is given by (\ref{phiCircle}) so that  $\Omega$ is the circle centered at $(0.5,0.5)$. Examples of meshes used both by $\phi$-FEM and by the standard FEM are given in Fig.~\ref{fig:meshes dirichlet}. We want to solve (\ref{eq:heat}) with the manufactured solution  \[u=u_{ex} = \exp(x)\sin(2\pi y)\sin(t) \] and extrapolated boundary conditions \[{u}^g = {u}_{ex}(1+\phi).   \]

We are going to compare the convergence of the $\phi$-FEM (\ref{discreteHet}) with that of the standard FEM using $\mathbb{P}^1$ Lagrange polynomials in space and the implicit Euler scheme in time in both cases.  The $\phi$-FEM stabilization parameter is taken as $\sigma = 20$. The results are reported in Figs.~\ref{fig:heat} and \ref{fig:heat2}, for $\Delta t=h$ and $\Delta t=10h^2$, respectively. Once again, $\phi$-FEM converges faster than standard FEM. In the test considered here, the predominant source of error seems to be in the time discretization. In particular, we observe only $O(h)$ convergence in the $L^2$-norm in space in the regime $\Delta t=h$ on Fig.~\ref{fig:heat}. A cleaner 2nd order in time should be possible to achieve using the BDF2 marching scheme, but this remains out of the scope of the present paper.

\FloatBarrier

\section{Conclusions and perspectives}


$\phi$-FEM is a relative newcomer to the field of unfitted FE methods. Up to now, it was only applied to scalar 2nd order elliptic equations with pure Dirichlet or pure Neumann/Robin boundary conditions in \cite{phifem,phiFEM2}. The purpose of the present contribution is to demonstrate its applicability to more sophisticated settings including the linear elasticity with mixed boundary conditions and material properties jumping across the internal interfaces, elasticity with cracks, and the heat transfer. In all the cases considered here, the numerical tests confirm the optimal accuracy on manufactured smooth solutions. $\phi$-FEM is easily implementable in standard FEM packages (we have chosen FEniCS for the numerical illustration in this chapter). 
In particular, $\phi$-FEM uses classical finite element spaces and avoids the mesh generation and any non-trivial numerical integration.

Interestingly, our methods systematically outperform the standard FEM on comparable meshes. This can be attributed to a better representation of the boundary and of the solution near the boundary, as opposed to the approximation of the domain by a polyhedron/polygon in standard FEM. We recall that the computing times, reported in some of our tests with $\phi$-FEM and favourably compared with those of the standard FEM, only include  assembling of the matrices and the resolution of the linear systems. It would be interesting to add the mesh generation time to the comparison, which should be even more in favour of $\phi$-FEM (when efficiently implemented).

Admittedly, the test cases presented in this contribution do not comprise all the complexity of the real-life problems. We have restricted ourselves to simple geometries in 2D only. Even more importantly, we have tested the methods only on smooth solutions, which is not supposed to happen in practice in problems with cracks, for example. Taking accurately into account the singularity at the crack tip remains an important challenge for the future $\phi$-FEM developments. A relatively easily implementable approach would be to combine $\phi$-FEM with a local mesh refinement by  quadtree/octree structures near the crack tip (front). We emphasize that such a refinement should be necessary only in the vicinity of the front, since the discontinuous solution along the crack should be efficiently approximated by $\phi$-FEM on a reasonably coarse unfitted mesh.  

The mathematical analysis of the schemes presented in this paper is in progress. We also plan to adapt $\phi$-FEM to fluid-structure simulations starting by the creeping flow around of a Newtonian fluid (Stokes equations) in the presence of rigid particles. 

\bibliographystyle{plain}

\bibliography{biblio}

\end{document}